\input ./style/arxiv-general.cfg
\documentclass[aos,MSNbibl,seceqn,dvips]{arximspdf}
\makeatletter
   \@ifpackageloaded{graphicx}{}{\usepackage{graphicx}}
\makeatother
\usepackage{dcolumn}

\doi{10.1214/15-AOS1327}
\volume{43}
\issue{4}
\pubyear{2015}
\firstpage{1831}
\lastpage{1864}
\docsubty{FLA}

\makeatletter
\newcolumntype{d}[1]{D{.}{.}{#1}}

\newcommand{\rrvert}{\vert}

\newcommand{\llvert}{\vert}

\newproclaim{defi}{Definition}
\newtheorem{theorem}{Theorem}
\newproclaim{remark}{Remark}
\newproclaim{exam}[prop]{Example}
\newtheorem{lemma}{Lemma}
\newproclaim{assumption}{Assumption}

\makeatother

\begin{document}
\begin{frontmatter}

\title{Jump activity estimation for pure-jump semimartingales via
self-normalized statistics\thanksref{T1}}
\runtitle{Jump activity estimation via self-normalized statistics}
\thankstext{T1}{Supported in part by NSF Grant SES-0957330.}

\begin{aug}
\author[A]{\fnms{Viktor}~\snm{Todorov}\corref{}\ead[label=e1]{v-todorov@northwestern.edu}}
\runauthor{V. Todorov}
\affiliation{Northwestern University}
\address[A]{Department of Finance\\
Northwestern University\\
Evanston, Illinois 60208-2001\\
USA\\
\printead{e1}}
\end{aug}

%
\received{\smonth{10} \syear{2014}}
%
\revised{\smonth{12} \syear{2014}}

%
\begin{abstract}
 We derive a nonparametric estimator of the jump-activity
index $\beta$ of a ``locally-stable'' pure-jump It\^o semimartingale
from discrete observations of the process on a fixed time interval with
mesh of the observation grid shrinking to zero. The estimator is based
on the empirical characteristic function of the increments of the
process scaled by local power variations formed from blocks of
increments spanning shrinking time intervals preceding the increments
to be scaled. The scaling serves two purposes: (1) it controls for the
time variation in the jump compensator around zero, and (2) it ensures
self-normalization, that is, that the limit of the characteristic
function-based estimator converges to a nondegenerate limit which
depends only on $\beta$. The proposed estimator leads to nontrivial
efficiency gains over existing estimators based on power variations. In
the L\'{e}vy case, the asymptotic variance decreases multiple times for
higher values of $\beta$. The limiting asymptotic variance of the
proposed estimator, unlike that of the existing power variation based
estimators, is constant. This leads to further efficiency gains in the
case when the characteristics of the semimartingale are stochastic.
Finally, in the limiting case of $\beta=2$, which corresponds to
jump-diffusion, our estimator of $\beta$ can achieve a faster rate than
existing estimators.
\end{abstract}

%

\pdfkeywords{62F12, 62M05, 60H10, 60J75, Central limit theorem,
high-frequency data, Ito semimartingale, jumps, jump activity index,
stochastic volatility, power variation}
%
\end{frontmatter}

\section{Introduction}
In this paper we are interested in estimating the jump activity index
of a process defined on a filtered probability space $(\Omega,\mathcal
{F},\break (\mathcal{F}_t)_{t\geq0},\mathbb{P})$ and given by
%
\begin{equation}
\label{eq:X} dX_t = \alpha_{t}\,dt+\sigma_{t-}\,dL_t+dY_t,
\end{equation}
when $L$ is locally stable pure-jump L\'{e}vy process (i.e., a pure-jump
L\'{e}vy process whose L\'{e}vy measure around zero behaves like that
of a
stable process) and $Y$ is a pure-jump process which is ``dominated''
at high-frequencies by $L$ in a sense which is made precise below; see
Assumption~\ref{assA}. All formal conditions for $X$ are given in Section~\ref
{sec:setting}. The jump activity index of $X$ on a given fixed time
interval is the infimum of the set of powers $p$ for which the sum of
$p$th absolute moments of the jumps is finite. Provided $\sigma$ does
not vanish on the interval and has c\`{a}dl\`{a}g paths, the jump activity
index of $X$ coincides with the Blumenthal--Getoor index of the driving
L\'{e}vy process $L$ (recall $Y$ is dominated by $L$ at high frequencies).
The dominant role of $L$ at high frequencies, together with its
stable-like L\'{e}vy measure around zero, manifests into the following
limiting behavior at high frequencies:
%
\begin{equation}\qquad
\label{eq:ls} h^{-1/\beta}(X_{t+sh}-X_t) \stackrel{
\mathcal{L}} {\longrightarrow } \sigma_{t}\times(S_{t+s}-S_t)\qquad
\mbox{as $h\rightarrow0$ and $s\in[0,1]$},
\end{equation}
for every $t$ and where $S$ is $\beta$-stable process, with the
convergence being for the Skorokhod topology. Equation~(\ref{eq:ls}) holds when
$\beta>1$ which is the case we consider in this paper. (When $\beta<1$
the drift will be the ``dominant'' component at high-frequencies, and
some of our results can be extended to this case as well.) We study
estimation of $\beta$ from discrete equidistant observations of $X$ on
a fixed time interval with mesh of the observation grid shrinking to zero.

Estimation of the jump activity index has received a lot of attention
recently. \cite{NR} consider estimation from low-frequency observations
in the setting of L\'{e}vy processes. \cite{B11} and \cite{BP13} consider
estimation from low-frequency data in the setting of time-changed L\'{e}vy
processes with an independent time-change process. \cite{B10}~consider
estimation from low-frequency and options data. \cite{B11_b} and \cite
{BP12} consider estimation from low frequency data in certain
stochastic volatility models. \cite{Woerner,Woerner03,Woerner07}
propose estimation from high-frequency data using power variations in a
pure-jump setting. \cite{SJ07} and \cite{JKLM12} consider estimation in
high-frequency setting when the underlying process can contain a
continuous martingale via truncated power variations. \cite{TT09}
propose estimation of the jump activity index in pure-jump setting via
power variations with adaptively chosen optimal power. \cite{T13}
extend \cite{TT09} via power variations of differenced increments which
provide further robustness and efficiency gains. \cite{JKL11} consider
jump activity estimation from noisy high-frequency data.

The estimation of $\beta$ from high-frequency data, thus far, makes use
of the dependence of the scaling factor of the high-frequency
increments in (\ref{eq:ls}) on $\beta$. For example, consider the power
variation
%
\begin{eqnarray}
\label{eq:pv} V(p,\Delta_n)& =& \sum_{i=1}^n\bigl|
\Delta_i^nX\bigr|^p, \qquad\Delta_i^n
X= X_{
{i}/{n}}-X_{{(i-1)}/{n}},
\nonumber
\\[-8pt]
\\[-8pt]
\nonumber
 \Delta_n &=& \frac{1}{n},\qquad
p>0.
\end{eqnarray}
Under certain technical conditions, (\ref{eq:ls}) implies
%
\begin{eqnarray*}
\Delta_n^{1-p/\beta}V(p,\Delta_n)& \stackrel{\mathbb
{P}} {\longrightarrow }& \mu\int_0^1|
\sigma_s|^p\,ds, \\
(2\Delta_n)^{1-p/\beta}V(p,2
\Delta _n) &\stackrel{\mathbb{P}} {\longrightarrow}& \mu\int
_0^1|\sigma _s|^p\,ds,
\end{eqnarray*}
where $\mu$ is some constant. An estimate of $\beta$ then can be simply
formed as a nonlinear function of the ratio $\frac{V(p,\Delta
_n)}{V(p,2\Delta_n)}$. This makes inference for $\beta$ possible
despite the unknown process $\sigma$.

The limit result in (\ref{eq:ls}), however, contains much more
information about $\beta$ than previously used in estimation. In
particular, (\ref{eq:ls}) implies that over a short interval of time
the increments of $X$, conditional on $\sigma$ at the beginning of the
interval, are approximately i.i.d. stable random variables. In this
paper we propose a new estimator of $\beta$ that utilizes this
additional information in (\ref{eq:ls}) and leads to significant
efficiency gains over existent estimators based on high-frequency data.

The key obstacle in utilizing the result in (\ref{eq:ls}) in inference
for $\beta$ is the fact that the process $\sigma$ is unknown and
time-varying. The idea of our method is to form a local estimator of
$\sigma$ using a block of high-frequency increments with asymptotically
shrinking time span via a localized version of (\ref{eq:pv}). We then
divide the high-frequency increments of $X$ by the local estimator of
$\sigma$. The division achieves ``self-normalization'' in the following
sense. First, the scale factor for the local estimator of $\sigma$ and
the high-frequency increment of $X$ are the same, and hence by taking
the ratio, they cancel. Second, both the high-frequency increment of
$X$ and the local estimator of $\sigma$ are approximately proportional
to the value of $\sigma$ at the beginning of the high-frequency
interval, and hence taking their ratio cancels the effect of the
unknown $\sigma$. The resulting scaled high-frequency increments are
approximately i.i.d. stable random variables, and we make inference for
$\beta$ via an analogue of the empirical characteristic function
approach, which has been used in various other contexts; see, for
example, \cite{CDH}.

After removing an asymptotic bias, the limit behavior of the empirical
characteristic function of the scaled high-frequency increments is
determined by two correlated normal random variables. One of them is
due to the limiting behavior of the empirical characteristic function
of the high-frequency increments scaled by the limit of the local power
variation. The other is due to the error in estimating the local scale
by the local realized power variation. Importantly, because of the
``self-normalization,'' the $\mathcal{F}$-conditional asymptotic
variance of the empirical characteristic function of the scaled
high-frequency increments is not random but rather a constant that
depends only on $\beta$ and the power $p$. This makes feasible
inference very easy.

When comparing the new estimator with existing ones based on the power
variation, we find nontrivial efficiency gains. There are two reasons
for the efficiency gains. First, as we noted above, our estimator makes
full use of the limiting result in (\ref{eq:ls}) and not just the
dependence of the scale of the high-frequency increments on $\beta$,
which is the case for existing ones. Second, by locally removing the
effect of the time-varying $\sigma$, we make the inference as if
$\sigma
$ is constant; that is, the limit variance is the same, regardless of
whether $X$ is L\'{e}vy or not. By contrast, the estimator based on the
ratios of power variations is asymptotically mixed normal with
$\mathcal
{F}$-conditional variance\vspace*{-2pt} of the form $K(p,\beta)\frac{\int_0^1|\sigma
_s|^{2p}\,ds}{ (\int_0^1|\sigma_s|^{p}\,ds )^2}$, for some
constant $K(p,\beta)$, and we note that $\frac{\int_0^1|\sigma
_s|^{2p}\,ds}{ (\int_0^1|\sigma_s|^{p}\,ds )^2}\geq1$ with
equality whenever the process $|\sigma|$ is almost everywhere constant
on the interval $[0,1]$. That is, the presence of time-varying $\sigma$
decreases the precision of the power-variation based estimator of
$\beta$.

The efficiency gains of our estimator are bigger for higher values of
$\beta$. In the limit case of $\beta=2$, which corresponds to $L$ being
a Brownian motion, we show that our estimator can achieve a faster rate
of convergence than the standard $\sqrt{n}$ rate for existing estimators.

The rest of the paper is organized as follows. In Section~\ref
{sec:setting} we introduce the setting. In Section~\ref{sec:constr} we
construct our statistic, and in Section~\ref{sec:limit} we derive its
limit behavior. In Section~\ref{sec:ja} we build on the developed limit
theory and construct new estimators of the jump activity and derive
their limit behavior. This section also shows the efficiency gains of
the proposed jump activity estimators over existing ones. Section~\ref
{sec:jd} deals with the limiting case of jump-diffusion. Sections~\ref
{sec:mc} and \ref{sec:emp} contain a Monte Carlo study and an empirical
application, respectively. Proofs are in Section~\ref{sec:proof}.

\section{Setting and assumptions}\label{sec:setting}
We start with introducing the setting and stating the assumptions that
we need for the results in the paper. We first recall that a L\'{e}vy
process $L$ with the characteristic triplet $(b,c,\nu)$, with respect
to truncation function $\kappa$ (Definition II.2.3 in \cite{JS}), is a
process with a characteristic function given by
%
\begin{eqnarray}
\label{eq:cf} \mathbb{E} \bigl(e^{iuL_t} \bigr) =\exp
\biggl[itub-tcu^2/2+t\int_{\mathbb
{R}}
\bigl(e^{iux}-1-iu\kappa(x) \bigr)\nu(dx) \biggr],
\nonumber
\\[-8pt]
\\[-8pt]
 \eqntext{t\geq0.}
\end{eqnarray}
In what follows we will always assume for simplicity that $\kappa
(-x)=-\kappa(x)$. Our assumption for the driving L\'{e}vy process in
(\ref{eq:X}) as well as the ``residual'' jump component $Y$ is given in
Assumption \ref{assA}.

\renewcommand{\theassumption}{\Alph{assumption}}
\begin{assumption}\label{assA}
$L$ in (\ref{eq:X}) is a
L\'{e}vy
process with characteristic triplet $(0,0,\nu)$ for $\nu$ a L\'{e}vy
measure with density given by
%
\begin{equation}
\label{eq:levy} \nu(x) = \frac{A}{|x|^{1+\beta}}+\nu'(x),\qquad \beta\in(0,2),
\end{equation}
where $A>0$ and $\nu'(x)$ is such that there exists $x_0>0$
with $|\nu'(x)|\leq C/|x|^{1+\beta'}$ for $|x|\leq x_0$ and some
$\beta
'<\beta$.

$Y$ is an It\^o semimartingale with the characteristic
triplet (\cite{JS}, Definition II.2.6) $ (\int_0^t\int_{\mathbb{R}}\kappa(x)\nu_s^Y(dx)\,ds,0,dt\otimes\nu_t^Y(dx)
)$ when
$\beta'<1$ and $ (0,0,dt\otimes\nu_t^Y(dx) )$ otherwise, with
$\int_{\mathbb{R}}(|x|^{\beta'+\iota}\wedge1)\nu_t^Y(dx)$ being
locally bounded and predictable, for some arbitrarily small $\iota>0$.
\end{assumption}

Assumption~\ref{assA} formalizes the sense in which $Y$ is dominated at high
frequencies by $L$: the activity index of $Y$ is below that of $L$. We
also stress that $Y$ and $L$ can have dependence. Therefore, as shown
in \cite{TT12}, we can accommodate in our setup time-changed L\'{e}vy
models, with absolute continuous time-change process, that have been
extensively used in applied work. Finally, we note that (\ref{eq:levy})
restricts only the behavior of $\nu$ around zero, and $\nu'$ is a
signed measure. Therefore many parametric jump specifications outside
of the stable process are satisfied by Assumption~\ref{assA} (e.g., the tempered
stable process). We next state our assumption for the dynamics of
$\alpha$ and $\sigma$.

\begin{assumption}\label{assB}
The processes $\alpha$ and
$\sigma$ are It\^o semimartingales of the form
%
\begin{eqnarray}
\label{ass:b} %
\alpha_t &=& \alpha_0+\int
_0^tb_s^{\alpha}\,ds+\int
_0^t\int_{E}\kappa
\bigl(\delta^{\alpha}(s,x)\bigr)\widetilde{\underline{\mu}}(ds,dx)\nonumber\\
&&{} +\int
_{E}\kappa'\bigl(\delta^{\alpha}(s,x)
\bigr)\underline{\mu}(ds,dx),
\nonumber
\\[-8pt]
\\[-8pt]
\nonumber
\sigma_t&=&\sigma_0+\int_0^tb_s^{\sigma}\,ds+
\int_0^t\int_{E}\kappa
\bigl(\delta ^{\sigma}(s,x)\bigr)\widetilde{\underline{\mu}}(ds,dx)\\
&&{} +\int
_{E}\kappa'\bigl(\delta^{\sigma}(s,x)
\bigr)\underline{\mu}(ds,dx), \nonumber
\end{eqnarray}
where $\kappa'(x) = x-\kappa(x)$, and:
\begin{longlist}[(a)]
\item[(a)] $|\sigma_t|^{-1}$ and $|\sigma_{t-}|^{-1}$ are
strictly positive;
\item[(b)] $\underline{\mu}$ is Poisson measure on $\mathbb
{R}_+\times E$, having arbitrary dependence with the jump measure of
$L$, with compensator $dt\otimes\lambda(dx)$ for some $\sigma$-finite
measures $\lambda$ on $E$;
\item[(c)] $\delta^{\alpha}(t,x)$ and $\delta^{\sigma}(t,x)$
are predictable, left-continuous with right limits in $t$ with $|\delta
^{\alpha}(t,x)|+|\delta^{\sigma}(t,x)|\leq\gamma_k(x)$ for all
$t\leq
T_k$, where $\gamma_k(x)$ is a deterministic function on $\mathbb{R}$
with $\int_{\mathbb{R}}(|\gamma_k(x)|^{r+\iota}\wedge1)\lambda
(dx)<\infty$ for arbitrarily small $\iota>0$ and some $0\leq r\leq
\beta
$, and $T_k$ is a sequence of stopping times increasing to $+\infty$;
\item[(d)] $b^{\alpha}$ and $b^{\sigma}$ are It\^o
semimartingales having dynamics as in (\ref{ass:b}) with coefficients
satisfying the analogues of conditions (b) and (c) above.
\end{longlist}
\end{assumption}

We note that $\underline{\mu}$ does not need to coincide with the jump
measure of $L$, and hence it allows for dependence between the
processes $\alpha$, $\sigma$ and $L$. This is of particular relevance
for financial applications. For example, Assumption~\ref{assB} is satisfied by
the COGARCH model of \cite{COGARCH} in which the jumps in $\sigma$ are
proportional to the squared jumps in $X$. More generally, Assumption~\ref{assB}
is satisfied if, for example, $(X,\alpha,\sigma)$ is modeled via a
L\'{e}vy-driven SDE, with each of the elements of the driving L\'{e}vy process
satisfying Assumption~\ref{assA}.

\section{Construction of the self-normalized statistics}\label{sec:constr}
We continue next with the construction of our statistics. The
estimation in the paper is based on observations of $X$ at the
equidistant grid times $0,\frac{1}{n},\ldots,1$ with $n\rightarrow
\infty$,
and we denote $\Delta_n = \frac{1}{n}$. To minimize the effect of the
drift in our statistics, we follow \cite{T13} and work with the first
difference of the increments, $\Delta_i^nX-\Delta_{i-1}^nX$, where
$\Delta_i^nX = X_{{i}/{n}}-X_{{(i-1)}/{n}}$ for $i=1,\ldots
,n$. The
above difference of increments is purged from the drift in the L\'{e}vy
case, and in the general case the drift has a smaller asymptotic effect
on it. For each $\Delta_i^nX-\Delta_{i-1}^nX$, we need a local power
variation estimate for the scale. It is constructed from a block of
$k_n$ high-frequency increments, for some $1<k_n<n-2$, as follows:
%
\begin{equation}
\label{eq:lpv} V_i^n(p) = \frac{1}{k_n}\sum
_{j=i-k_n-1}^{i-2}\bigl|\Delta_j^nX-
\Delta _{j-1}^nX\bigr|^p,\qquad  i=k_n+3,
\ldots,n.
\end{equation}
Block-based local estimators of volatility have been also used in other
contexts in a high-frequency setting, for example, in \cite{JR} and
\cite{TT14}. The empirical characteristic function of the scaled
differenced increments is given by
%
\begin{equation}
\label{eq:ecf} \widehat{\mathcal{L}}^n(p,u) = \frac{1}{n-k_n-2}\sum
_{i=k_n+3}^n\cos \biggl(u\frac{\Delta_i^nX-\Delta_{i-1}^nX}{ (V_i^n(p))^{1/p} }
\biggr),\qquad u\in\mathbb{R}_+.
\end{equation}
We proceed with some notation needed for the limiting theory of
$\widehat{\mathcal{L}}^n(p,u)$. Let $S_1$, $S_2$ and $S_3$ be random
variables corresponding to the values of three independent L\'{e}vy
processes at time 1, each of which with the characteristic triplet
$(0,0,\nu)$, for any truncation function $\kappa$ and where $\nu$ has
the density $\frac{A}{|x|^{1+\beta}}$. Then we denote
$\mu_{p,\beta} = (\mathbb{E}|S_1-S_2|^p)^{\beta/p}$, which does not
depend on $\kappa$, and we further use the shorthand notation $\mathbb
{E} ( e^{iu(S_1-S_2)} ) = e^{-A_{\beta}u^{\beta}}$ for any
$u>0$ with $A_{\beta}$ being a (known) function of $A$ and $\beta$.
Using Example~25.10 in \cite{Sato} and references therein, we have
%
\begin{equation}
\label{eq:C} C_{p,\beta} = \frac{A_{\beta}}{\mu_{p,\beta}}= \biggl[\frac
{2^{p}\Gamma
 ({(1+p)}/{2} )\Gamma (1-{p}/{\beta}
)}{\sqrt
{\pi}\Gamma (1-{p}/{2} )}
\biggr]^{-\beta/p},
\end{equation}
which depends only on $p$ and $\beta$ but not on the scale parameter of
the stable random variables $S_1$ and $S_2$.
With this notation, we set
%
\begin{equation}
\label{eq:limit} \mathcal{L}(p,u,\beta) = e^{-C_{p,\beta}u^{\beta}},\qquad u\in\mathbb{R}_+,
\end{equation}
which will be the limit in probability of $\widehat{\mathcal
{L}}^n(p,u)$. We finish with some more notation needed to describe the
asymptotic variance of $\widehat{\mathcal{L}}^n(p,u)$. First, we denote
for some $u\in\mathbb{R}_+$,
%
\begin{eqnarray}
\xi_1(p,u,\beta) &=& \biggl( %
\cos \biggl(\frac{u(S_1-S_2)}{\mu_{p,\beta
}^{1/\beta
}} \biggr) - \mathcal{L}(p,u,
\beta),  \frac{|S_1-S_2|^p}{\mu
_{p,\beta
}^{p/\beta}}-1
 \biggr)',
\nonumber
\\[-8pt]
\\[-8pt]
\nonumber
\xi_2(p,u,\beta) &=& \biggl( %
 \cos
\biggl(\frac{u(S_2-S_3)}{\mu_{p,\beta
}^{1/\beta
}} \biggr) - \mathcal{L}(p,u,\beta),
\frac{|S_2-S_3|^p}{\mu
_{p,\beta
}^{p/\beta}}-1
 \biggr)'.
\end{eqnarray}
We then set for $u,v\in\mathbb{R}_+$
%
\begin{equation}
\Xi_i(p,u,v,\beta) = \mathbb{E} \bigl(\xi_1(p,u,\beta)
\xi '_{1+i}(p,v,\beta) \bigr),\qquad i=0,1
\end{equation}
and
%
\begin{eqnarray}
G(p,u,\beta)& =& \frac{\beta}{p}e^{-C_{p,\beta}u^{\beta}}C_{p,\beta
}u^{\beta},
\nonumber
\\[-8pt]
\\[-8pt]
\nonumber
H(p,u,\beta)& =& G(p,u,\beta) \biggl( \frac{\beta
}{p}C_{p,\beta
}u^{\beta}
-\frac{\beta}{p}-1 \biggr).
\end{eqnarray}
%

\section{Limit theory for \texorpdfstring{$\widehat{\mathcal{L}}^n(p,u)$}{widehat{mathcal{L}}n(p,u)}}\label
{sec:limit} We start with convergence in probability.

\begin{theorem}\label{thm:cp}
Assume $X$ satisfies Assumptions \ref{assA} and \ref{assB} for some $\beta\in(1,2)$ and
$\beta'<\beta$. Let $k_n$ be a deterministic sequence satisfying
$k_n\asymp n^{\varpi}$ for some $\varpi\in(0,1)$. Then, for
$0<p<\beta
$, we have
%
\begin{equation}
\label{cp_1} \widehat{\mathcal{L}}^n(p,u) \stackrel{\mathbb{P}} {
\longrightarrow } \mathcal{L}(p,u,\beta) \qquad\mbox{as $n\rightarrow\infty$},
\end{equation}
locally uniformly in $u\in\mathbb{R}_+$.
\end{theorem}

We note that we restrict $\beta>1$; that is, we focus on the infinite
variation case. The above theorem will continue to hold for $\beta\leq
1$, but for the subsequent results about the limiting distribution of
$\widehat{\mathcal{L}}^n(p,u)$, we will need quite stringent additional
restrictions in the case $\beta\leq1$. We do not pursue this here. The
other conditions for the convergence in probability result are weak.
The requirements for $\alpha$ and $\sigma$ for Theorem~\ref{thm:cp} to
hold are actually much weaker than what is assumed in Assumption~\ref{assB}, but
for simplicity of exposition we keep Assumption~\ref{assB} throughout. We note
that for consistency, we have a lot of flexibility about the block size
$k_n$: (1)~$k_n\rightarrow\infty$ so that we consistently estimate the
scale via $V_i^n(p)$ and (2) $k_n/n\rightarrow0$ so that the span of
the block is asymptotically shrinking to zero, and therefore no bias is
generated due to the time variation of $\sigma$. In the case when $X$
is a L\'{e}vy process, the second condition is obviously not needed.

To derive a central limit theorem (c.l.t.) for $\widehat{\mathcal
{L}}^n(p,u)$, we will need to restrict the choice of $k_n$ more. We
will assume $k_n/\sqrt{n}\rightarrow0$, so that biases due to the time
variation in $\sigma$, which are hard to feasibly estimate, are
negligible. For such a choice of $k_n$, however, an asymptotic bias due
to the sampling error of $V_i^n(p)$ appears, and for stating a c.l.t., we need to
consider the following bias-corrected estimator:
%
\begin{eqnarray}
\label{eq:ecf_debias} \qquad
\widehat{\mathcal{L}}^{n}(p,u,
\beta)' &= &\widehat{\mathcal {L}}^{n}(p,u)
\nonumber
\\[-8pt]
\\[-8pt]
\nonumber
&&{} - \frac{1}{k_n}\frac{1}{2}H(p,u,\beta) \bigl( \Xi
_0^{(2,2)}(p,u,u,\beta) +2\Xi_1^{(2,2)}(p,u,u,
\beta) \bigr). %
\end{eqnarray}
We state the c.l.t. for $\widehat{\mathcal{L}}^{n}(p,u,\beta)' $ in the
next theorem.

\begin{theorem}\label{thm:clt}
Assume $X$ satisfies Assumptions \ref{assA} and \ref{assB} with $\beta\in(1,2)$ and
$\beta
'<\frac{\beta}{2}$, and that the power $p$ and block size $k_n$ satisfy
%
\begin{eqnarray}
\label{clt_1} \frac{\beta\beta'}{2(\beta-\beta')}&\vee&\frac{\beta
-1}{2}<p<\frac
{\beta}{2},
\\
\label{clt_2} k_n\asymp n^{\varpi},\qquad \frac{p}{\beta}&\vee&
\frac{1}{3}<\varpi <\frac{1}{2}.
\end{eqnarray}
Then, as $n\rightarrow\infty$, we have
%
\begin{equation}
\label{clt_3} \sqrt{n} \bigl( \widehat{\mathcal{L}}^n(p,u,
\beta)' - \mathcal {L}(p,u,\beta) \bigr) \stackrel{\mathcal{L}} {
\longrightarrow } Z_1(u)+G(p,u,\beta)Z_2(u),
\end{equation}
locally uniformly in $u\in\mathbb{R}_+$. $Z_1(u)$ and $Z_2(u)$ are two
Gaussian processes with the following covariance structure:
%
\begin{equation}
\label{clt_4} \mathbb{E} \bigl(\mathbf{Z}(u)\mathbf{Z}(v) \bigr) = \Xi
_0(p,u,v,\beta )+2\Xi_1(p,u,v,\beta),\qquad u,v\in
\mathbb{R}_+,
\end{equation}
where $\mathbf{Z}(u) =  (Z_1(u), Z_2(u) )'$.

Let $\widehat{\beta}$ be an estimator of $\beta$ with $\widehat
{\beta
}-\beta= o_p(k_n\sqrt{\Delta_n})$ as $n\rightarrow\infty$. Then
%
\begin{equation}
\label{clt_5} \sqrt{n} \bigl( \widehat{\mathcal{L}}^n(p,u,\widehat{
\beta})' - \widehat {\mathcal{L}}^n(p,u,
\beta)' \bigr) \stackrel{\mathbb {P}} {\longrightarrow} 0,
\end{equation}
locally uniformly in $u\in\mathbb{R}_+$.
\end{theorem}

The conditions for the power $p$ in (\ref{clt_1}) are exactly the same
as in \cite{T13} for the analysis of the realized power variation, and
they are relatively weak. For example, the condition $p>\frac{\beta
-1}{2}$ will be always satisfied as soon as we pick power slightly
above $\frac{1}{2}$. Moreover, this condition is not needed in the case
when $X$ is a L\'{e}vy process. Further, the condition in (\ref
{clt_2}) for
$k_n$ shows that we have more flexibility for the choice of $k_n$
whenever $p$ is not very close to its upper bound of $\beta/2$.

Due to the self-normalization in the construction of our statistic, the
limiting distribution in (\ref{clt_3}) is Gaussian and not mixed
Gaussian, which is the case for most limit results in high-frequency
asymptotics (and in particular for the power variation based estimator
of $\beta$); see \cite{V12} for another exception. This is very
convenient as the estimation of the asymptotic variance is
straightforward. The bias correction in (\ref{eq:ecf_debias}) is
infeasible, as it depends on $\beta$. However, (\ref{clt_5}) shows that
a feasible version of the debiasing would work provided\vspace*{1pt} the initial
estimator of $\beta$ is $o_p(k_n\sqrt{\Delta_n})$. When one estimates
$\beta$ using $\widehat{\mathcal{L}}^n(p,u)$, with explicit estimators
provided in the next section, $\widehat{\beta}-\beta$ will be
$O_p(1/k_n)$. Hence, such a preliminary estimate of $\beta$ will
satisfy the required rate condition in Theorem~\ref{thm:clt}.

\section{Jump activity estimation}\label{sec:ja} We now use the limit
theory developed above to form estimators of $\beta$. The simplest one
is based on $\widehat{\mathcal{L}}^n(p,u)$ and is given by
%
\begin{equation}
\label{beta:fs} \widehat{\beta}^{fs}(p,u,v) = \frac{\log (-\log(\widehat
{\mathcal
{L}}^n(p,u)) )-\log (-\log(\widehat{\mathcal
{L}}^n(p,v))
)}{\log(u/v)},
\end{equation}
for $u,v\in\mathbb{R}_+$ with $u\neq v$. Because of the asymptotic bias
in $\widehat{\mathcal{L}}^n(p,u)$, $\widehat{\beta
}^{fs}(p,\break u,v)-\beta$
will be only $O_p(1/k_n)$, with $p$ and $k_n$ satisfying (\ref
{clt_1})--(\ref{clt_2}). An explicit estimate of $\beta$ using feasible
debiasing is given by
%
\begin{equation}\qquad
\label{beta:ts} \widehat{\beta}(p,u,v) = \frac{\log (-\log(\widehat{\mathcal
{L}}^n(p,u,\widehat{\beta}^{fs})') )-\log (-\log
(\widehat
{\mathcal{L}}^n(p,v,\widehat{\beta}^{fs})') )}{\log(u/v)},
\end{equation}
for some $u,v\in\mathbb{R}_+$ with $u\neq v$, and where
$\widehat{\beta}^{fs}$ is a suitable initial estimator of $\beta$ [like
the one in (\ref{beta:fs})]. While convenient, the above estimators
have two potential drawbacks. One, we do not take into account the
information about $\beta$ in the constant $C_{p,\beta}$. This is
because in the asymptotic limit of the above estimators, $C_{p,\beta}$
gets canceled. Second, $u$ and $v$ are chosen arbitrarily, and one can
include more moment conditions for the estimation of $\beta$ using
$\widehat{\mathcal{L}}^n(p,u,\widehat{\beta}^{fs})'$. In the next
theorem we provide a general estimator of $\beta$ which overcomes these
drawbacks of the explicit estimators above.

\begin{theorem}\label{thm:gmm}
Assume $X$ satisfies Assumptions \ref{assA} and \ref{assB} with $\beta\in(1,2)$ and
$\beta
'<\beta/2$, and that the conditions in (\ref{clt_1}) and (\ref{clt_2})
hold. Suppose $\widehat{\beta}^{fs}$ is a consistent estimator of
$\beta
$ with $\widehat{\beta}^{fs}-\beta= o_p(k_n\sqrt{\Delta_n})$. Denote
with $\widehat{\mathbf{u}}_l$ and $\widehat{\mathbf{u}}_h$ two
sequences of $K\times1$-dimensional vectors, for some finite $K\geq
1$, satisfying $\widehat{\mathbf{u}}_l \stackrel{\mathbb
{P}}{\longrightarrow} \mathbf{u}_l$ and $\widehat{\mathbf
{u}}_h \stackrel{\mathbb{P}}{\longrightarrow} \mathbf{u}_h$ as
$n\rightarrow\infty$, for some $\mathbf{u}_l, \mathbf{u}_h\in
\mathbb
{R}_+^K$ with $u_l^i<u_h^i$, $u_l^j<u_h^j$ and $(u_l^i, u_h^i)\cap
(u_l^j, u_h^j)= \varnothing$ for every $i,j = 1,\ldots,K$ with $i\neq j$
where $u_l^i$ and $u_h^i$ denote the $i$th element of the vectors
$\mathbf{u}_l$ and $\mathbf{u}_h$, respectively. Set further the
shorthand $\mathbf{u}=[\mathbf{u}_l; \mathbf{u}_h]$ and $\widehat
{\mathbf{u}}=[\widehat{\mathbf{u}}_l; \widehat{\mathbf{u}}_h]$.

Let $\mathbf{W}(p,\mathbf{u},\beta)$ be $K\times K$ matrix with $(i,j)$
element given by
%
\begin{eqnarray}
\label{gmm_1} \mathbf{W}(p,\mathbf{u},\beta)_{i,j} &=& \int
_{u_l^i}^{u_h^i}\int_{u_l^j}^{u_h^j}w(p,u,v,
\beta)\,du\,dv,
\\
\label{gmm_2} %
w(p,u,v,\beta) &= &\frac{1}{ \mathcal{L}(p,u,\beta)\mathcal
{L}(p,v,\beta)
}\pmatrix{ 1
\cr
G(p,u,\beta)}'\nonumber\\
&&{}\times\overline{
\Xi}(p,u,v,\beta) \pmatrix{ 1
\cr
G(p,v,\beta)},\nonumber
\end{eqnarray}
where $\overline{\Xi}(p,u,v,\beta) = \Xi_0(p,u,v,\beta)+2\Xi
_1(p,u,v,\beta) $.

Define the $K\times1$ vector $\widehat{\mathbf{m}}(p,\widehat
{\mathbf
{u}},\widehat{\beta}^{fs},\mathbf{u},\beta)$ by
%
\begin{equation}\quad
\label{gmm_3} \widehat{\mathbf{m}}\bigl(p,\widehat{\mathbf{u}},\widehat{\beta
}^{fs},\mathbf {u},\beta\bigr)_i = \int
_{\widehat{u}_l^i}^{\widehat{u}_h^i} \bigl(\log \bigl(\widehat{
\mathcal{L}}^n\bigl(p,u,\widehat{\beta}^{fs}
\bigr)'\bigr) - \log \bigl(\mathcal {L}(p,u,\beta)\bigr) \bigr)\,du,
\end{equation}
for $i=1,\ldots,K$, and set
%
\begin{eqnarray}
\label{eq:beta_gmm} &&\widehat{\beta}(p,\mathbf{u})
\nonumber
\\[-8pt]
\\[-8pt]
\nonumber
&&\qquad= \mathop{\operatorname{argmin}}_{\beta\in(1,2)}
\widehat{\mathbf{m}}\bigl(p,\widehat{\mathbf {u}},\widehat{\beta
}^{fs},\mathbf{u},\beta\bigr)'\mathbf{W}^{-1}
\bigl(p,\widehat{\mathbf {u}},\widehat {\beta}^{fs}\bigr) \widehat{
\mathbf{m}}\bigl(p,\widehat{\mathbf{u}},\widehat {\beta }^{fs},
\mathbf{u},\beta\bigr).
\end{eqnarray}
Finally define the $K\times1$ vector $\mathbf{M}(p,\mathbf{u},\beta
)$ by
%
\begin{equation}
\label{gmm_5} \mathbf{M}(p,\mathbf{u},\beta)_i = \int
_{u_l^i}^{u_h^i}\nabla _{\beta
}\log\bigl(
\mathcal{L}(p,u,\beta)\bigr)\,du,\qquad i=1,\ldots,K.
\end{equation}

Then for $\beta\in(1,2)$, $p\in (\frac{\beta\beta'}{2(\beta
-\beta
')},\frac{\beta}{2} )$ and $\beta'<\beta/2$, we have
%
\begin{equation}
\label{gmm_6} %
\sqrt{n} \bigl( \widehat{\beta}(p,\mathbf{u}) - \beta
\bigr) \stackrel{\mathcal{L}} {\longrightarrow} \sqrt{\mathbf{M}(p,\mathbf {u},
\beta)' \mathbf{W}^{-1}(p,\mathbf{u},\beta)\mathbf {M}(p,
\mathbf {u},\beta)}\times\mathcal{N}, %
\end{equation}
for $n\rightarrow\infty$ with $\mathcal{N}$ being standard normal
random variable.

A consistent estimator for the asymptotic variance of $\widehat{\beta
}(p,\mathbf{u})$ is given by
%
\begin{equation}
\label{gmm_7} \mathbf{M}(p,\widehat{\mathbf{u}},\widehat{\beta})'
\mathbf {W}^{-1}(p,\widehat{\mathbf{u}},\widehat{\beta}) \mathbf {M}(p,
\widehat {\mathbf{u}},\widehat{\beta}),
\end{equation}
where $\mathbf{M}(p,\widehat{\mathbf{u}},\widehat{\beta})$ is defined
as $\mathbf{M}(p,\mathbf{u},\beta)$ with $\mathbf{u}$ and $\beta$
replaced by $\widehat{\mathbf{u}}$ and $\widehat{\beta}$.
\end{theorem}

Theorem~\ref{thm:gmm} allows us to adaptively choose the range of $u$
over which to match $\widehat{\mathcal{L}}^n(p,u,\widehat{\beta
}^{fs})'$ with its limit. This is convenient because the limiting
variance of $\widehat{\mathcal{L}}^n(p,u,\widehat{\beta}^{fs})'$
depends on $\beta$. For this reason also the weight function in (\ref
{gmm_1}) optimally weighs the moment conditions in the estimation. We
discuss the practical issues regarding the construction of $\widehat
{\mathbf{m}}(p,\widehat{\mathbf{u}},\widehat{\beta}^{fs},\mathbf
{u},\beta)$ in Section~\ref{sec:mc}.

\begin{figure}

\includegraphics{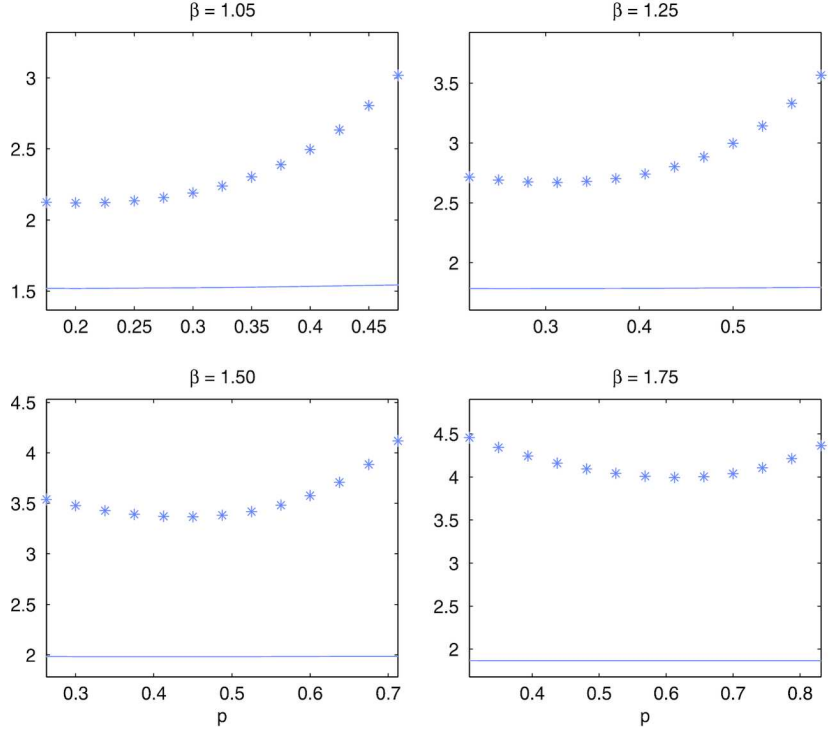}

\caption{Asymptotic standard deviation of jump activity
estimators. The straight line corresponds to the asymptotic standard
deviation of the characteristic function based estimator defined in
(\protect\ref{eq:beta_gmm}) and the $*$ line to the power variation based
estimator of \cite{T13} given in (\protect\ref{eq:beta_pv}) (when
$\sigma$ is
constant). For each cases of $\beta$, the power $p$ ranges in the
interval $p\in (\frac{7}{40},\frac{19}{40} )\beta$. For the
estimator in (\protect\ref{eq:beta_gmm}), the vector $\mathbf{u}_l =
[0.1:0.05:5]$ and $\mathbf{u}_h = [0.15:0.05:5.05]$.}
\label{fig:beta_ase}
\end{figure}

We now illustrate the efficiency gains provided by the new method over
existing power variation based estimators of $\beta$. The power
variation estimator based on the differenced increments is given by
(see \cite{T13})
%
\begin{equation}
\label{eq:beta_pv} \widetilde{\beta}(p) = \frac{p\log(2)}{\log [\widetilde
{V}_2^n(p)/\widetilde{V}_1^n(p) ]}1_{\{\widetilde{V}_1^n(p)\neq
\widetilde{V}_2^n(p)\}},
\end{equation}
where
%
\begin{eqnarray}
\widetilde{V}_1^n(p) &=& \sum
_{i=2}^n\bigl|\Delta_i^nX-\Delta
_{i-1}^nX\bigr|^p,
\nonumber
\\[-8pt]
\\[-8pt]
\nonumber
 \widetilde{V}_2^n(p)
&= &\sum_{i=4}^n\bigl|\Delta_i^nX-
\Delta _{i-1}^nX+\Delta_{i-2}^nX-
\Delta_{i-3}^nX\bigr|^p.
\end{eqnarray}
On Figure~\ref{fig:beta_ase}, we plot the limiting standard deviation
of the estimators in (\ref{eq:beta_gmm}) and (\ref{eq:beta_pv}) for
different values of $\beta$. [The estimator in (\ref{eq:beta_pv}) is
derived under exactly the same assumptions for $X$ as our estimator
here.] The asymptotic standard deviation of $\widetilde{\beta}(p)$ is
computed from \cite{T13}. $\widehat{\beta}(p,u)$ is far less sensitive
to the choice of $p$ than $\widetilde{\beta}(p)$, with lower powers
yielding marginally more efficient $\widehat{\beta}(p,u)$. The news
estimator $\widehat{\beta}(p,u)$ provides nontrivial efficiency gains
irrespective of the values of $p$ and $\beta$. The gains are bigger for
high values of the jump activity. For example, for $\beta=1.75$,
$\widehat{\beta}(p,u)$ is around two times more efficient (in terms of
asymptotic standard deviation) than $\widetilde{\beta}(p)$.

\section{The limiting case of jump-diffusion}\label{sec:jd}
So far our analysis has been for the pure-jump case of $\beta\in(1,2)$.
We now look at the limiting case of $\beta=2$, which corresponds to $L$
in (\ref{eq:X}) being a Brownian motion. In this case the asymptotic
behavior of the high-frequency increments in (\ref{eq:ls}) holds with
$S$ being a Brownian motion. Thus deciding $\beta=2$ versus $\beta<2$
amounts to testing pure-jump versus jump-diffusion specification for
$X$. It turns out that when $\beta=2$, our estimation method can lead
to a faster rate of convergence than the $\sqrt{n}$ rate we have seen
for the case $\beta\in(1,2)$. This is unlike the power-variation based
estimation methods for which the rate of convergence is $\sqrt{n}$,
both for $\beta=2$ and $\beta<2$; see, for example,~\cite{TT09}.

The faster rate of convergence in the case $\beta=2$ can be achieved by
letting the argument $u$ of the empirical characteristic function
$\widehat{\mathcal{L}}(p,u)$ drift toward zero as $n\rightarrow
\infty$.
In this case, $\frac{-\log(\widehat{\mathcal
{L}}(p,u_n,2)')}{C_{p,2}u_n^2}$ and $\frac{-\log(\widehat{\mathcal
{L}}(p,\rho u_n,2)')}{C_{p,2}\rho^2u_n^2}$, for some $\rho>0$, are
asymptotically perfectly correlated, and their difference converges at
a faster rate. We note that this does not work in the pure-jump case of
$\beta<2$. To state the formal result we first introduce some notation.
For $S_1$, $S_2$ and $S_3$ being independent standard normal random
variables, we denote
%
\begin{eqnarray}
\widetilde{\xi}_1(p) &=& \biggl( %
 \frac{|S_1-S_2|^4}{\mu_{p,2}^{2}}- \frac{12}{\mu
_{p,2}^{2}},  \frac{|S_1-S_2|^p}{\mu_{p,2}^{p/2}}-1
 \biggr)',
\nonumber
\\[-8pt]
\\[-8pt]
\nonumber
\widetilde{\xi}_2(p) &= &\biggl( %
\frac{|S_2-S_3|^4}{\mu_{p,2}^{2}}- \frac{12}{\mu
_{p,2}^{2}},  \frac{|S_2-S_3|^p}{\mu_{p,2}^{p/2}}-1
 \biggr)', %
\end{eqnarray}
and then set $\widetilde{\Xi}_i(p) = \mathbb{E} (\widetilde
{\xi
}_1(p)\widetilde{\xi}'_{1+i}(p) )$ for $i=0,1$. The difference
from the analogous expression for the case $\beta<2$ is in the first
terms of $\widetilde{\xi}_1(p)$ and $\widetilde{\xi}_2(p)$. Note that
the expression for the bias-correction remains exactly the same as it
involves only the variance and covariance of the second elements of
$\widetilde{\xi}_1(p)$ and $\widetilde{\xi}_2(p)$, which remain the
same as their pure-jump counterparts.

\begin{theorem}\label{thm:cont}
Suppose $X$ has dynamics given by (\ref{eq:X}) with $L$ being a
Brownian motion, $Y$ satisfying the corresponding condition for it in
Assumption~\ref{assA} and $\alpha$ and $\sigma$ satisfying Assumption~\ref{assB} for some
$r<2$. Suppose $p<1$, $k_n\sqrt{\Delta_n}\rightarrow0$ and
$u_n\rightarrow0$, and further
%
\begin{eqnarray}
\label{cont_1} \frac{ \Delta_n^{  ({p}/{\beta'}-{p}/{2}
)\wedge
{(p+1)}/{(r\vee1+1)} -\iota}\vee k_n^{- ({1}/{p}\wedge
{3}/{2} )+\iota} \vee(k_n\Delta_n)^{1-\iota}}{u_n^6\sqrt
{\Delta
_n}}&\rightarrow&0,
\nonumber
\\[-8pt]
\\[-8pt]
\nonumber
 \frac{(k_n\Delta_n)^{{1}/{r}\wedge
{(2-p)}/{2}-{1}/{2}}}{u_n^6}&\rightarrow& 0.
\end{eqnarray}
Then for some $\rho>0$
%
\begin{equation}
\label{cont_2} \widehat{\beta}^{fs}(p,u_n,\rho
u_n) - 2 = O_p\bigl(k_n^{-1}u_n^2
\bigr).
\end{equation}
Further, if for some initial estimator $\widehat{\beta}^{fs} -2 =
o_p(k_nu_n^{2}\sqrt{\Delta_n})$, then
%
\begin{equation}\qquad
\label{cont_3} \frac{\sqrt{n}}{u_n^2(1-\rho^2)} \bigl(\widehat{\beta}(p,u_n,\rho
u_n)-2 \bigr) \stackrel{\mathcal{L}} {\longrightarrow} -
\frac
{1}{\log
(\rho)} \biggl(\frac{1}{24C_{p,2}}Z_1-\frac{2}{p}C_{p,2}Z_2
\biggr),
\end{equation}
where $Z_1$ and $Z_2$ are two zero-mean normal random variables with
covariance given by
$\widetilde{\Xi}_0(p)+2\widetilde{\Xi}_1(p)$.

When $X$ is a L\'{e}vy process, the requirement for $k_n$ and $u_n$
reduces to
%
\begin{equation}
\label{cont_4} u_n\rightarrow0,\qquad \frac{ \Delta_n^{  ({p}/{\beta
'}\wedge
1-{p}/{2} ) -\iota}\vee k_n^{- ({1}/{p}\wedge
{3}/{2} )+\iota} }{u_n^6\sqrt{\Delta_n}}\rightarrow0.
\end{equation}
\end{theorem}

The rate of convergence of the estimator for $\beta$ is now $\sqrt
{n}u_n^{-2}$ and is faster than the one in Theortem \ref{thm:gmm}, when
$u_n$ converges to zero. The latter is determined by the restriction in
(\ref{cont_1}), which in turn is governed by the presence of the
``residual'' term $Y$, the variation in $\sigma$ and the sampling
variation in measuring the scale via $V_i^n(p)$. For the condition to
be satisfied we need $p\in(1/2,1)$ and $\beta'<1$; that is, the jumps
in $X$ are of finite variation; for testing the null hypothesis of
presence of diffusion when the process can contain infinite variation
jumps under the null, see the recent work of \cite{JKL14}. Without any
prior knowledge on $\beta'$ and $r$, we can set $k_n$ according to
(\ref{clt_2}), with $\beta=2$, and then set $u_n\asymp\log(n)^{-1}$.
The requirement on $u_n$ can be further relaxed when $X$ is a L\'{e}vy
process as evident from (\ref{cont_4}). Finally, we can draw a parallel
between our finding for faster rate of convergence of the estimator of
$\beta$ when $\beta=2$ with the result in \cite{DM73,DM83} for faster
rate of convergence for the maximum likelihood estimator of the
stability index of i.i.d. $\beta$-stable random variables when $\beta=2$.

\section{Monte Carlo}\label{sec:mc}
We test the performance of the proposed method for jump activity
estimation on simulated data from the following model
%
\begin{equation}
dX_t = \sigma_{t-}\,dL_t,\qquad d\sigma_t
= -0.03 \sigma_t\,dt+dZ_t,
\end{equation}
where $L$ and $Z$ are two L\'{e}vy processes independent of each other with
L\'{e}vy densities given by $\nu_L(x) = e^{-\lambda|x|} (\frac
{A_0}{|x|^{1+\beta}}+\frac{A_1}{|x|^{1+{\beta}/{3}}} )$ and
$\nu_Z(x) = \break 0.0293\frac{e^{-3x}}{x^{1.5}}1_{\{x>0\}}$, respectively.
$\sigma$ is a L\'{e}vy-driven Ornstein--Uhlenbeck process with a tempered
stable driving L\'{e}vy subordinator. The parameters governing the dynamics
of $\sigma$ imply $\mathbb{E}(\sigma_t) = 1$ and half-life of shock in
$\sigma$ of around one month (when unit of time is a day). $L$ is a
mixture of tempered stable processes with the parameter $\beta$
coinciding with the jump activity index of $X$. We fix $\lambda=
0.25$, and consider four cases for~$\beta$. In each of the cases we set
$A_0$ and $A_1$ so that $A_0\int_{\mathbb{R}}|x|^{1-\beta
}e^{-\lambda
|x|}\,dx = 1$ and $A_1\int_{\mathbb{R}}|x|^{1-{\beta
}/{3}}e^{-\lambda
|x|}\,dx = 0.2$. The four cases are: (1) $\beta= 1.05$ and $A_0 =
0.1299$, $A_1 = 0.0113$; (2) $\beta= 1.25$ and $A _0= 0.1443$, $A_1 =
0.0125$; (3) $\beta= 1.50$ and $A_0 = 0.1410$, $A_1 = 0.0141$ and (4)
$\beta= 1.75$ and $A_0 = 0.0975$, $A_1 = 0.0158$.

In the Monte Carlo we set $T = 10$ and $n = 100$ which corresponds
approximately to two weeks of 5-minute return data in a typical
financial setting. We further set $k_n = 50$ and $p=0.51$. The initial
estimator to construct the moments and the optimal weight matrix is
simply $\widehat{\beta}^{fs}(p,u,v)$ with $u = 0.1$ and $v=1.1$. If
$p\geq\widehat{\beta}^{fs}(p,u,v)/2$, then we reduce the power to $p =
\widehat{\beta}^{fs}(p,u,v)/4$. Based on the initial beta estimator, we
estimate the values of $u$ for which $\mathcal{L}(p,u,\beta)=0.95$ and
$\mathcal{L}(p,u,\beta)=0.25$, and then split this interval in five
equidistant regions which are used in constructing the moment vector in
(\ref{gmm_2}).

Regarding the number of moment conditions, $K$, in the construction of
our estimator, we should keep in mind the following. Larger $K$ helps
improve efficiency of the estimator as our equal weighting of the
characteristic function within each moment condition is suboptimal.
However, the feasible estimate of the optimal weight matrix is unstable
in small samples when $K$ is large. (This is similar to ``curse of
dimensionality'' problems occurring in related contexts; see, e.g.,
\cite{FLK} and \cite{K}.) Moreover, since the characteristic function
is smooth, one typically does not need many moment conditions to gain
efficiency. For example, we also experimented in the Monte Carlo with
ten moment conditions (by splitting the region of $u$ into ten
equidistant regions). The performance of the estimator based on the ten
moment conditions was very similar to the one based on the five moment
conditions whose performance we summarize below.

The results from the Monte Carlo are reported in Table~\ref{tb:mc}. For
comparison, we also report results for $\widetilde{\beta}(p)$ where $p$
is set to the level which minimizes the corresponding asymptotic
standard deviation in Figure~\ref{fig:beta_ase}. We notice satisfactory
finite sample performance of $\widehat{\beta}(p,\mathbf{u})$. In all
cases for $\beta$, $\widehat{\beta}(p,\mathbf{u})$ contains relatively
small upward biases. These biases, however, are well below those of
$\widetilde{\beta}(p)$. We note that the finite sample bias of
$\widehat
{\beta}(p,\mathbf{u})$ can be significantly reduced if, similar to
$\widetilde{\beta}(p)$, one uses an adaptive choice of power in the
range $(\beta/4,\beta/3)$. The superiority of $\widehat{\beta
}(p,\mathbf
{u})$ holds also in terms of precision in estimating $\beta$, with
inter-quantile ranges of $\widehat{\beta}(p,\mathbf{u})$ typically well
below those of $\widetilde{\beta}(p)$.

\begin{table}
\caption{Monte Carlo results}\label{tb:mc}
\begin{tabular*}{\textwidth}{@{\extracolsep{\fill}}lcccccc@{}}
\hline
 & \multicolumn{3}{c}{$\bolds{\widehat{\beta}(p,\mathbf{u})}$} &
 \multicolumn{3}{c}{$\bolds{\widetilde{\beta}(p)}$}\\[-6pt]
  & \multicolumn{3}{c}{\hrulefill} &
 \multicolumn{3}{c@{}}{\hrulefill}\\
\textbf{Case} & \textbf{Median} & \textbf{IQR} & \textbf{MAD} &  \textbf{Median} & \textbf{IQR} & \textbf{MAD}\\
\hline
$\beta= 1.05$ & $1.0801$ & $0.0791$ & $0.0518$ & $1.1154$ &
$0.0925$ & $0.0792$\\
$\beta= 1.25$ & $1.3058$ & $0.0817$ & $0.0680$ & $1.3229$ &
$0.1158$ & $0.0932$\\
$\beta= 1.50$ & $1.5398$ & $0.0886$ & $0.0622$ & $1.5767$ &
$0.1405$ & $0.1072$\\
$\beta= 1.75$ & $1.7782$ & $0.0806$ & $0.0536$ & $1.8196$ &
$0.1704$ & $0.1183$\\
\hline
\end{tabular*}
\tabnotetext[]{}{\textit{Note}: IQR is the inter-quartile range, and MAD is the mean
absolute deviation around the true value. The power $p$ for $\widetilde
{\beta}(p)$ is set to the value which minimizes the corresponding
asymptotic standard deviation displayed in Figure~\ref{fig:beta_ase}.}
\end{table}

\section{Empirical application}\label{sec:emp}
\begin{figure}

\includegraphics{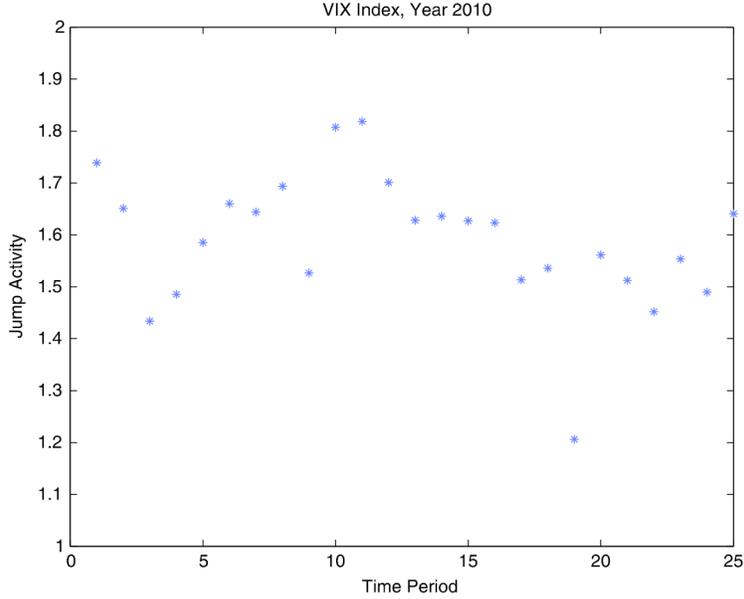}

\caption{Jump Activity for the VIX Index. Estimation is done
over periods of $10$ days in the year $2010$. In the estimation,
moments $p$ and $k_n$ are selected as in the Monte Carlo.}
\label{fig:emp}
\end{figure}

We now apply the developed inference procedures on high-frequency data
for the VIX index. The VIX index is a option-based measure for
volatility in the market (S\&P 500 index). It serves as a popular
indicator for investors' uncertainty, and it is used as the underlying
asset for many volatility-based derivative contracts traded in the
financial exchanges. Earlier work, consistent with parametric models
for volatility, has provided evidence that the VIX index is a pure-jump
It\^o semimartingale. Here, we estimate its jump activity index. The
estimation is based on $5$-minute sampled data during the trading hours
for the year $2010$. Like in the Monte Carlo, we split the year into
intervals of $10$ days (two weeks) and estimate the jump activity over
each of them. The moments, the power $p$ and the block size $k_n$, are
selected in the same way as in the Monte Carlo. Estimation results are
presented in Figure~\ref{fig:emp}. The estimated jump activity index
takes values around $1.6$. Overall, our results support a pure-jump
specification of the VIX index.
%

\section{Proofs}\label{sec:proof}
In the proofs we use the shorthand notation $\mathbb{E}_i^n(\cdot)
\equiv\mathbb{E}(\cdot|\mathcal{F}_{i\Delta_n})$ and $\mathbb
{P}_i^n(\cdot) \equiv\mathbb{P}(\cdot|\mathcal{F}_{i\Delta_n})$. We
also denote with $K$ a positive constant that does not depend on $n$
and $u$ and might change from line to line in the inequalities that
follow. When we want to highlight that the constant depends only on
some parameters $a$ and $b$, we write $K_{a,b}$.

\subsection{Decompositions and additional notation}
In what follows it is convenient to extend appropriately the
probability space and then decompose the driving L\'{e}vy process $L$
as follows:
%
\begin{equation}
\label{proof:jd} L_t+\widehat{S}_t = S_t+
\widetilde{S}_t,
\end{equation}
where $S$, $\widehat{S}$ and $\widetilde{S}$ are pure-jump L\'{e}vy
processes with the first two characteristics zero [with respect to the
truncation function $\kappa(\cdot)$] and L\'{e}vy densities $\frac
{A}{|x|^{1+\beta}}$, $2|\nu'(x)|1_{\{\nu'(x)<0\}}$ and $|\nu'(x)|$,
respectively. We denote the associated counting jump measures with $\mu
$, $\mu_1$ and $\mu_2$. (Note that there can be dependence between
$\mu
$, $\mu_1$ and $\mu_2$.)

$S$ is $\beta$-stable process, and $\widehat{S}$ and $\widetilde{S}$
are ``residual'' components whose effect on our statistic, as will be
shown, is negligible (under suitable conditions). The proof of the
decomposition in (\ref{proof:jd}) as well as the explicit construction
of $S$, $\widehat{S}$ and $\widetilde{S}$ can be found in Section~1 of
the supplementary Appendix of \cite{TT12}.

We now introduce some additional notation that will be used throughout
the proofs. We denote for $i=k_n+3,\ldots,n$,
%
\begin{eqnarray*}
\widehat{V}_i^n(p) &=& \frac{1}{k_n}\sum
_{j=i-k_n-1}^{i-2}|\sigma _{(j-2)\Delta_n-}|^p\bigl|
\Delta_j^nS-\Delta_{j-1}^nS\bigr|^p,
\\
\overline{V}_i^n(p) &= &\frac{1}{k_n}\sum
_{j=i-k_n-1}^{i-2}\frac
{|\Delta
_j^nS-\Delta_{j-1}^nS|^p}{\mu_{p,\beta}^{p/\beta}},
\\
\dot{V}_i^n(p) &=& \sum
_{j=i-k_n-1}^{i-2} \biggl\{\frac{[(i-j-4)\vee
0+1_{\{j<i-3\}}]}{k_n} \bigl(|
\sigma_{j\Delta_n-}|^p - |\sigma _{(j-2)\Delta_n-}|^p
\bigr)
\\
&&\hspace*{103pt}{} + \frac{ (|\sigma_{(j-1)\Delta_n-}|^p-|\sigma_{(j-2)\Delta
_n-}|^p)1_{\{
j<i-2\}} }{k_n} \biggr\}\\
&&\hspace*{38pt}{}\times \bigl|\Delta_j^nS-
\Delta_{j-1}^nS\bigr|^p,
\\
|\overline{\sigma}|_i^p& =& \frac{1}{k_n}\sum
_{j=i-k_n-1}^{i-2}|\sigma _{(j-2)\Delta_n-}|^p.
\end{eqnarray*}
We further denote the function
%
\[
f_{i,u}(x) = \exp \biggl( -\frac{C_{p,\beta}u^{\beta}|\sigma
_{(i-2)\Delta
_n-}|^{\beta}}{x^{\beta/p} } \biggr),
\]
and direct computation yields
%
\[
\cases{\displaystyle %
f_{i,u}'(x) =
\frac{\beta}{p}f_{i,u}(x)\frac{
C_{p,\beta
}u^{\beta}|\sigma_{(i-2)\Delta_n-}|^{\beta} }{x^{\beta/p+1} },
\vspace*{2pt}\cr
\displaystyle f_{i,u}^{\prime\prime}(x) = f_i(u,x) \biggl(
\frac{\beta}{p}\frac{
C_{p,\beta
}u^{\beta}|\sigma_{(i-2)\Delta_n-}|^{\beta} }{x^{\beta/p+1} } \biggr)^2 \vspace*{2pt}\cr
\displaystyle \hspace*{43pt}{}- f_i(u,x)
\frac{\beta}{p} \biggl( \frac{\beta}{p}+1 \biggr)\frac{
C_{p,\beta}u^{\beta}|\sigma_{(i-2)\Delta_n-}|^{\beta} }{ x^{\beta/p+2}
}.}
\]
We note
%
\begin{equation}
\label{decomp_1} \sup_{x\in\mathbb{R}_+}\bigl|f_{i,u}(x) +
f_{i,u}'(x) + f_{i,u}^{\prime\prime}(x) +
f_{i,u}^{\prime\prime\prime}(x)\bigr|<K_u,
\end{equation}
where the positive constant $K_u$ depends only on $u$ and is finite as
soon as $u$ is bounded away from zero.

With this notation, we make the following decomposition for any $u\in
\mathbb{R}_+$:
%
\[
\widehat{\mathcal{L}}^n(p,u) - \mathcal{L}(p,u,\beta) =
\frac
{1}{n-k_n-2} \Biggl[\widehat{Z}_1^n(u)+
\widehat{Z}_2^n(u)+\sum_{j=1}^4R_j^n(u)
\Biggr],
\]
where $\widehat{Z}_j^n(u) = \sum_{i=k_n+3}^nz_i^j(u)$ for $j=1,2$ with
%
\begin{eqnarray*}
z_i^1(u)& =& \cos \biggl(u\frac{\sigma_{(i-2)\Delta_n-}(\Delta
_i^nS-\Delta
_{i-1}^nS)}{ (V_i^n(p))^{1/p} } \biggr) -
\exp \biggl(-\frac{A_{\beta
}u^{\beta}|\sigma_{(i-2)\Delta_n-}|^{\beta} }{\Delta
_n^{-1}(V_i^n(p))^{\beta/p}} \biggr),
\\
z_i^2(u)&= & \exp \biggl(-\frac{C_{p,\beta}u^{\beta}|\sigma
_{(i-2)\Delta
_n-}|^{\beta} }{\Delta_n^{-1} ( |\overline{\sigma
}|_i^p\overline
{V}_i^n(p) )^{\beta/p}} \biggr) -
\exp \biggl( -\frac
{C_{p,\beta
}u^{\beta}|\sigma_{(i-2)\Delta_n-}|^{\beta}}{ (|\overline
{\sigma
}|_i^p  )^{\beta/p}} \biggr),
\end{eqnarray*}
and $R_j^n(u) = \sum_{i=k_n+3}^nr_i^j(u)$ for $j=1,2,3,4$ with
%
\begin{eqnarray*}
r_i^1(u)& =&\cos \biggl(u\frac{\Delta_i^nX-\Delta_{i-1}^nX}{
(V_i^n(p))^{1/p} } \biggr) - \cos
\biggl(u\frac{\sigma_{(i-2)\Delta
_n-}(\Delta_i^nS-\Delta_{i-1}^nS)}{ (V_i^n(p))^{1/p} } \biggr),
\\
r_i^2(u) &=& \exp \biggl(-\frac{A_{\beta}u^{\beta}|\sigma
_{(i-2)\Delta
_n-}|^{\beta} }{\Delta_n^{-1}(V_i^n(p))^{\beta/p}} \biggr) -
\exp \biggl(-\frac{A_{\beta}u^{\beta}|\sigma_{(i-2)\Delta_n-}|^{\beta}
}{\Delta
_n^{-1}(\widehat{V}_i^n(p))^{\beta/p}} \biggr),
\\
r_i^3(u)& =& \exp \biggl(-\frac{A_{\beta}u^{\beta}|\sigma
_{(i-2)\Delta
_n-}|^{\beta} }{\Delta_n^{-1}(\widehat{V}_i^n(p))^{\beta/p}} \biggr) -
\exp \biggl(- \frac{C_{p,\beta}u^{\beta}|\sigma_{(i-2)\Delta
_n-}|^{\beta
}}{\Delta_n^{-1} ( |\overline{\sigma}|_i^p\overline
{V}_i^n(p)
)^{\beta/p} } \biggr),
\\
r_i^4(u) &=& \exp \biggl( -\frac{C_{p,\beta}u^{\beta}|\sigma
_{(i-2)\Delta
_n-}|^{\beta}}{ (|\overline{\sigma}|_i^p  )^{\beta/p}} \biggr) -
\exp \bigl(-C_{p,\beta}u^{\beta} \bigr).
\end{eqnarray*}
We finally introduce the following: $\overline{Z}_1^n(u) = \sum_{i=k_n+3}^n\overline{z}_i^1(u)$,
$\overline{Z}_2^{(a,n)}(u) = \break  \sum_{i=k_n+3}^n\overline{z}_i^{(a,2)}(u)$ and $\overline{Z}_2^{(b,n)}(u) =
\sum_{i=k_n+3}^n\overline{z}_i^{(b,2)}(u)$ where
%
\begin{eqnarray*}
\overline{z}_i^{1}(u) &=& \cos \bigl(u\Delta_n^{-1/\beta}
\mu _{p,\beta
}^{-1/\beta}\bigl(\Delta_i^nS-
\Delta_{i-1}^nS\bigr) \bigr) - \mathcal {L}(p,u,\beta),
\\
\overline{z}_i^{(a,2)}(u)& =&G(p,u,\beta) \bigl(
\Delta_n^{-p/\beta
}\overline{V}_i^n(p) -1
\bigr),\\
 \overline{z}_i^{(b,2)}(u) &= &\tfrac
{1}{2}H(p,u,
\beta) \bigl( \Delta_n^{-p/\beta}\overline{V}_i^n(p)
-1 \bigr)^2.
\end{eqnarray*}
%

\subsection{Localization}
We prove results under the following strengthened version of Assumption~\ref{assB}:

\renewcommand{\theassumption}{SB}
\begin{assumption}\label{assSB} We have Assumption \ref{assB} and in addition:
\begin{longlist}[(a)]
\item[(a)] the processes $|\sigma_t|$ and $|\sigma_t|^{-1}$
are uniformly bounded;
\item[(b)] the processes $b^{\alpha}$ and $b^{\sigma}$ are
uniformly bounded;
\item[(c)] $|\delta^{\alpha}(t,x)|+|\delta^{\sigma
}(t,x)|\leq
\gamma(x)$ for all $t$, where $\gamma(x)$ is a deterministic bounded
function on $\mathbb{R}$ with $\int_{\mathbb{R}}|\gamma
(x)|^{r+\iota
}\lambda(dx)<\infty$ for arbitrarily small $\iota>0$ and some $0\leq
r\leq\beta$;
\item[(d)] the coefficients in the It\^o semimartingale
representation of $b^{\alpha}$ and $b^{\sigma}$ satisfy the analogues
of conditions (b) and (c) above;
\item[(e)] the process $\int_{\mathbb{R}}(|x|^{\beta
'+\iota
}\wedge1)\nu_t^Y(dx)$ is bounded, and the jumps of $\widehat{S}$,
$\widetilde{S}$ and $Y$ are bounded.
\end{longlist}
\end{assumption}

Extending the results to the case of the more general Assumption~\ref{assB}
follows by standard localization arguments given in Section~4.4.1 of
\cite{JP}.

\subsection{Preliminary results}
The strategy of the proofs is to bound the terms $R_j^n(u)$ for
$j=1,2,3,4$ as well as $\widehat{Z}_1^n(u) - \overline{Z}_1^n(u)$ and
$\widehat{Z}_2^n(u) - \overline{Z}_2^{(a,n)}(u) - \overline
{Z}_2^{(b,n)}(u)$, and to derive the asymptotic limits of $\overline
{Z}_1^n(u)$, $\overline{Z}_2^{(a,n)}(u)$ and $\overline
{Z}_2^{(b,n)}(u)$. We do this in a sequence of lemmas starting with one
containing some preliminary bounds needed for the subsequent lemmas.

\begin{lemma}\label{lema:prelim-a}
Under Assumptions \ref{assA} and \ref{assSB} and $k_n\asymp n^{\varpi}$ for $\varpi\in
(0,1)$, we have for $0<p<\beta$, $\iota>0$ arbitrarily small and
$1\leq
x<\frac{\beta}{p}$ and $y\geq1$,
%
\begin{eqnarray}
\label{prelim_a-1} %
&&\Delta_n^{-p/\beta}
\mathbb{E}\bigl|V_i^n(p) - \widehat{V}_i^n(p)\bigr|
\leq K\alpha_n, %
\\
\nonumber
&&\alpha_n = \frac{\Delta_n^{(2-1/\beta)(1+(p-1/2)\wedge0-\iota
)}}{\sqrt
{k_n}}\vee\Delta_n^{{1}/{\beta}-\iota}
\vee\Delta _n^{
{p}/{\beta'}\wedge1-{p}/{\beta}-\iota}\\
&&\hspace*{24pt}{}\vee\Delta_n^{
{(p+1)}/{(\beta+1)}-\iota},\nonumber
\\
\label{prelim_a-2}
&&\mathbb{E}\bigl\llvert \Delta_n^{-p/\beta}
\widehat{V}_i^n(p) - \mu _{p,\beta
}^{p/\beta}|
\overline{\sigma}|_i^p\bigr\rrvert ^{x}+
\mathbb{E}\bigl\llvert \Delta_n^{-p/\beta}\overline{V}_i^n(p)
- 1\bigr\rrvert ^{x}
\nonumber
\\[-8pt]
\\[-8pt]
\nonumber
&&\qquad\leq K\cases{ %
k_n^{-x/2},&\quad $\mbox{if $\beta/p>2$},$
\vspace*{2pt}\cr
k_n^{1-x},&\quad $\mbox{if $\beta/p\leq2$},$}
\\
\label{prelim_a-3}
&&\bigl\llvert \mathbb{E}_{i-k_n-3}^n\bigl(|
\overline{\sigma}|_i^p - |\sigma _{(i-2)\Delta_n-}|^p
\bigr)\bigr\rrvert \leq Kk_n\Delta_n,
\\
\label{prelim_a-3-b}
&& \mathbb{E}_{i-k_n-3}^n\bigl\llvert |\overline{
\sigma}|_i^p - |\sigma _{(i-2)\Delta_n-}|^p
\bigr\rrvert ^y\leq K(k_n\Delta_n)^{{y}/{r}\wedge
1-\iota},
\\
\label{prelim_a-4}
&&\Delta_n^{-p/\beta}\bigl\llvert
\mathbb{E}_{i-k_n-3}^n \bigl(\widehat {V}_i^n(p)
- \mu_{p,\beta}^{p/\beta}|\overline{\sigma }|_i^p
\overline {V}_i^n(p) - \dot{V}_i^n(p)
\bigr) \bigr\rrvert \leq Kk_n\Delta_n,
\\
\label{prelim_a-4-b}
&& \Delta_n^{-xp/\beta}\mathbb{E}\bigl\llvert
\widehat{V}_i^n(p) - \mu _{p,\beta
}^{p/\beta}|
\overline{\sigma}|_i^p\overline{V}_i^n(p)
- \dot {V}_i^n(p)\bigr\rrvert ^{x}\leq
K(k_n\Delta_n)^{{x}/{r}\wedge
1-\iota},
\\
\label{prelim_a-4-c}
&& \Delta_n^{-xp/\beta}\mathbb{E}\bigl\llvert
\dot{V}_i^n(p)\bigr\rrvert ^{x}\leq K
\Delta_n^{{(\beta-xp)}/{\beta}\wedge{x}/{r}-\iota}.
\end{eqnarray}
\end{lemma}
\begin{pf}
We start with
(\ref{prelim_a-1}). We apply exactly the same decomposition and bounds
as for the term $A_3$ in Section~5.2.3 in \cite{T13} to get the result
in (\ref{prelim_a-1}).
We continue with (\ref{prelim_a-2}). Without loss of generality we
assume $k_n\geq2$, and we denote the two sets
%
\[
\label{proof_a_1} \cases{ %
\displaystyle J_i^{e}
= \biggl\{i-k_n-1+2k\dvtx k=0,\ldots,\biggl\lfloor \frac
{k_n-1}{2}
\biggr\rfloor \biggr\},
\vspace*{2pt}\cr
\displaystyle J_i^{o} = \biggl\{i-k_n-1+2k+1\dvtx k=0,
\ldots,\biggl\lfloor\frac{k_n-2}{2}\biggr\rfloor \biggr\}.}
\]
With this notation, we can decompose $\widehat{V}_i^{n}(p)$ into
%
\begin{eqnarray*}
\label{proof_a_2} \widehat{V}_i^{(e,n)}(p)& =& \frac{1}{k_n}
\sum_{j\in J_i^e}|\sigma _{(j-2)\Delta_n-}|^p\bigl|
\Delta_j^nS-\Delta_{j-1}^nS\bigr|^p,\\
\widehat {V}_i^{(o,n)}(p) &=& \widehat{V}_i^{n}(p)
- \widehat {V}_i^{(e,n)}(p).
\end{eqnarray*}
We further denote $|\overline{\sigma}|_{e,i}^p = 
\frac{1}{k_n}\sum_{j\in J_i^e}|\sigma_{(j-2)\Delta_n-}|^p$ and $|\overline{\sigma
}|_{o,i}^p =\break   \frac{1}{k_n}\sum_{j\in J_i^o}|\sigma_{(j-2)\Delta
_n-}|^p$. Using the triangular inequality, we then have
%
\begin{eqnarray*}
\label{proof_a_3} %
&&\bigl\llvert \Delta_n^{-p/\beta}
\widehat{V}_i^{n}(p) - \mu_{p,\beta
}^{p/\beta
}|
\overline{\sigma}|_i^p\bigr\rrvert
\\
&&\qquad \leq \bigl\llvert \Delta_n^{-p/\beta}\widehat{V}_i^{(e,n)}(p)
- \mu_{p,\beta
}^{p/\beta}|\overline{\sigma}|_{e,i}^p
\bigr\rrvert + \bigl\llvert \Delta_n^{-p/\beta}
\widehat{V}_i^{(o,n)}(p) - \mu_{p,\beta
}^{p/\beta}|
\overline{\sigma}|_{o,i}^p\bigr\rrvert .
\end{eqnarray*}
Now, since $\mathbb{E}_{j-2}^n|\Delta_j^nS-\Delta_{j-1}^nS|^p =
\Delta
_n^{p/\beta}\mu_{p,\beta}^{p/\beta}$, the sums $\Delta_n^{-p/\beta
}\widehat{V}_i^{(e,n)}(p) - \mu_{p,\beta}^{p/\beta}|\overline
{\sigma
}|_{e,i}^p$ and $\Delta_n^{-p/\beta}\widehat{V}_i^{(o,n)}(p) - \mu
_{p,\beta}^{p/\beta}|\overline{\sigma}|_{o,i}^p$ are discrete
martingales. From here, the result in (\ref{prelim_a-2}) for the case
$\beta/p\leq2$ follows by a direct application of the
Burkholder--Davis--Gundy inequality and the algebraic inequality
%
\begin{equation}
\label{proof_a_3-a} \biggl|\sum_i|a_i|\biggr|^p
\leq\sum_{i}|a_i|^p\qquad
\forall p\in(0,1] \mbox{ and any real-valued $\{a_{i}
\}_{i\geq1}$}.
\end{equation}
We are left with the case $\beta/p>2$. We only show the bound involving
the term $\widehat{V}_i^{(e,n)}(p)$, with the result for $\widehat
{V}_i^{(o,n)}(p)$ being shown analogously. We first denote $\Delta
_n^{-p/\beta}\widehat{V}_i^{(e,n)}(p) - \mu_{p,\beta}^{p/\beta
}|\overline{\sigma}|_{e,i}^p = \frac{1}{k_n}\sum_{j\in J_i^e}\zeta_j^n$
where $\zeta_j^n = \Delta_n^{-p/\beta}\times\break |\sigma_{(j-2)\Delta
_n-}|^p  (
|\Delta_j^nS-\Delta_{j-1}^nS|^p - \mu_{p,\beta}^{p/\beta}  )$.
Applying the Burkholder--Davis--Gundy inequality, we have
%
\[
\mathbb{E}\biggl| \sum_{j\in J_i^e}\zeta_j^n
\biggr|^x \leq K\mathbb{E} \biggl(\sum_{j\in J_i^e}
\bigl(\zeta_j^n\bigr)^2 \biggr)^{x/2}.
\]
If $x\leq2$, the result in (\ref{prelim_a-2}) then follows by Jensen's
inequality. If $x>2$, applying again Burkholder--Davis--Gundy, we have
%
\begin{eqnarray}
\label{proof_a_3-b} %
&&\mathbb{E} \biggl(\sum_{j\in J_i^e}
\bigl(\zeta_j^n\bigr)^2 \biggr)^{x/2}
\nonumber\\
&&\qquad \leq K\mathbb{E} \biggl(\sum_{j\in J_i^e}\bigl(\bigl(
\zeta_j^n\bigr)^2-\mathbb
{E}_{j-2}^n\bigl(\zeta_j^n
\bigr)^2\bigr) \biggr)^{x/2}+K\mathbb{E} \biggl(\sum
_{j\in
J_i^e}\mathbb{E}_{j-2}^n\bigl(
\zeta_j^n\bigr)^2 \biggr)^{x/2}
\\
&&\qquad \leq K\mathbb{E} \biggl(\sum_{j\in J_i^e}\bigl(\bigl(
\zeta_j^n\bigr)^2-\mathbb
{E}_{j-2}^n\bigl(\zeta_j^n
\bigr)^2\bigr)^2 \biggr)^{x/4}+Kk_n^{x/2},\nonumber
\end{eqnarray}
where we also made use of the fact that the $\beta$-stable random
variable has finite $p$th absolute moment as soon as $p\in(0,\beta)$.
If $x\leq4$, the result will then follow from an application of (\ref
{proof_a_3-a}). If $x>4$, then we repeat (\ref{proof_a_3-b}) with $x$
replaced by $x/2$ and $\zeta_j^n$ replaced $(\zeta_j^n)^2-\mathbb
{E}_{j-2}^n(\zeta_j^n)^2$. We continue in this way, applying $k = \sup
 \{i\dvtx2^i<x \}$ times (\ref{proof_a_3-b}) and then
(\ref{proof_a_3-a}). This shows (\ref{prelim_a-2}).

We continue with (\ref{prelim_a-3}) and (\ref{prelim_a-3-b}). We make
use of the following algebraic inequality:
%
\[
\label{proof_a_4} \bigl\llvert |a+b|^p-|a|^p - p
\operatorname{sign}\{a\}|a|^{p-1}b\bigr\rrvert \leq K_p|a|^{p-2}|b|^2,
\]
for any $a,b\in{R}$ with $a\neq0$, $0<p<1$ and $K_p$ that depends only
on $p$. Applying this inequality as well as the triangular inequality,
and using the fact that under Assumption \ref{assSB} the process $|\sigma|$ is
bounded from below, we have
%
\begin{eqnarray}
\label{proof_a_5} \bigl\llvert \mathbb{E}_s\bigl(|
\sigma_t|^p-|\sigma_s|^p\bigr)
\bigr\rrvert &\leq& K|t-s|,\qquad 0\leq s\leq t,
\\
\mathbb{E}_s\bigl\llvert |\sigma_t|^p-|
\sigma_s|^p\bigr\rrvert ^q&\leq& K\mathbb
{E}_s \bigl(|\sigma_t-\sigma_s|^q
\vee|\sigma_t-\sigma _s|^{2q} \bigr),
\nonumber
\\[-8pt]
\\[-8pt]
\eqntext{ 0\leq s
\leq t, q\geq1,}
\end{eqnarray}
with some constant $K$ that does not depend on $s$ and $t$. From here
(\ref{prelim_a-3}) follows. Application of Corollary~2.1.9 of \cite{JP}
further gives
%
\begin{equation}
\label{proof_a_6} \mathbb{E}_s|\sigma_t-
\sigma_s|^q\leq K|t-s|^{{q}/{r}\wedge
1-\iota},\qquad 0\leq s\leq t,
q\geq1,
\end{equation}
and applying this inequality with $q=y$ and $q=2y$, for $y$ the
constant in (\ref{prelim_a-3-b}), we have that result.

We proceed by showing the bounds in (\ref{prelim_a-4})--(\ref
{prelim_a-4-c}). We can decompose $|\overline{\sigma}|_i^p - |\sigma
_{(k-2)\Delta_n-}|^p = \sum_{j=1}^4a_k^j$ for $k=i-k_n-1,\ldots,i-2$ and
%
\[
\label{proof_a_7} \cases{ %
\displaystyle a_k^1
= \frac{1}{k_n}\sum_{j=k+3}^{i-2} \bigl(|
\sigma _{(j-2)\Delta_n-}|^p-|\sigma_{k\Delta_n-}|^p
\bigr),
\vspace*{2pt}\cr
\displaystyle a_k^2 = \frac
{(i-k-4)\vee0}{k_n} \bigl(|\sigma_{k\Delta_n-}|^p
- |\sigma _{(k-2)\Delta
_n-}|^p \bigr),
\vspace*{2pt}\cr
\displaystyle a_k^3 = \bigl(\bigl (|\sigma_{k\Delta_n-}|^p - |\sigma
_{(k-2)\Delta_n-}|^p\bigr)1_{\{k<i-3\}}\vspace*{2pt}\cr
\hspace*{28pt}{}+\bigl(|\sigma_{(k-1)\Delta
_n-}|^p-|\sigma
_{(k-2)\Delta_n-}|^p\bigr)1_{\{k<i-2\}} \bigr)/{k_n},
\vspace*{2pt}\cr
\displaystyle a_k^4 = \frac{1}{k_n}\sum
_{j=i-k_n-1}^k \bigl(|\sigma_{(j-2)\Delta
_n-}|^p-|
\sigma_{(k-2)\Delta_n-}|^p \bigr),}
\]
with $a_k^1$ being zero for $k\geq i-4$. Using the law of iterated
expectations and the bound in (\ref{prelim_a-3-b}), we have for
$k=i-k_n-1,\ldots,i-2$,
%
\begin{equation}
\label{proof_a_8} \Delta_n^{-xp/\beta}\mathbb{E} \bigl(\bigl\llvert
a_k^1+a_k^4\bigr\rrvert \bigl|
\Delta _k^nS-\Delta_{k-1}^nS\bigr|^p
\bigr)^x\leq K(k_n\Delta_n)^{
{x}/{r}\wedge1-\iota}.
\end{equation}
Using the H\"{o}lder inequality, the bound in (\ref{proof_a_5}), as well
as the fact that a stable random variable has finite absolute moments
for powers less than $\beta$, we have for $k=i-k_n-1,\ldots,i-2$,
%
\begin{eqnarray}
\label{proof_a_9} &&\Delta_n^{-xp/\beta}\mathbb{E} \bigl(\bigl\llvert
a_k^2+a_k^3\bigr\rrvert\bigl |
\Delta _k^nS-\Delta_{k-1}^nS\bigr|^p
\bigr)^x
\nonumber
\\[-8pt]
\\[-8pt]
\nonumber
&&\qquad\leq K\Delta_n^{ ({(\beta x/r)}/{(\beta-xp)}\wedge1 ){(\beta-xp)}/{\beta} -\iota}.
\end{eqnarray}
Combining (\ref{proof_a_8}) and (\ref{proof_a_9}), we get the results
in (\ref{prelim_a-4-b}) and (\ref{prelim_a-4-c}).

Further, using (\ref{proof_a_5}), we get for $k=i-k_n-1,\ldots,i-2$,
%
\begin{equation}
\label{proof_a_10} \Delta_n^{-p/\beta}\bigl\llvert
\mathbb{E}_{i-k_n-3} ^n \bigl(\bigl(a_k^1+a_k^4
\bigr)\bigl|\Delta_k^nS-\Delta_{k-1}^nS\bigr|^p
\bigr)\bigr\rrvert \leq Kk_n\Delta_n.
\end{equation}
From here we get the result in (\ref{prelim_a-4}).
\end{pf}

\begin{lemma}\label{lema:prelim-b}
Under Assumptions \ref{assA} and \ref{assSB} and $k_n\asymp n^{\varpi}$ for $\varpi\in
(0,1)$, we have for $0<p<\beta$, $\iota>0$ arbitrarily small and every
$0<a<b<\infty$,
%
\begin{eqnarray}\qquad\hspace*{4pt}
\label{prelim_b-1}
&& \frac{1}{n-k_n-2}\mathbb{E} \Bigl(\sup_{u\in[a,b]}\bigl|R_1^n(u)\bigr|
\Bigr)
\leq K_{a,b} \bigl(\alpha_n\vee k_n^{- ({\beta}/{(2p)}\wedge{(\beta-p)}/{p} )+\iota}
\bigr),
\\
\label{prelim_b-2}
&&\frac{1}{n-k_n-2}\mathbb{E} \Bigl(\sup_{u\in[a,b]}\bigl|R_2^n(u)\bigr|
\Bigr)\leq K_{a,b}\alpha_n,
\\
\label{prelim_b-3} %
&&\frac{1}{n-k_n-2}\mathbb{E} \Bigl(\sup
_{u\in[a,b]}\bigl|R_3^n(u)\bigr| \Bigr)
\nonumber
\\[-8pt]
\\[-8pt]
\nonumber
&&\qquad \leq\cases{ %
\displaystyle K_{a,b}
\bigl((k_n\Delta_n)^{1-\iota}\vee
k_n^{-1/2}(k_n\Delta_n)^{{1}/{r}\wedge{(\beta-p)}/{\beta
}-\iota
}
\bigr),&\quad $\mbox{if $\beta/p>2$,}$
\vspace*{2pt}\cr
\displaystyle K_{a,b} \bigl((k_n\Delta _n)^{
{1}/{r}\wedge1-\iota}
\vee\Delta_n^{{(\beta-p)}/{\beta}-\iota
} \bigr),&\quad $\mbox{if $\beta/p\leq2$,}$}
\\
\label{prelim_b-4}
&& \frac{1}{n-k_n-2}\mathbb{E} \Bigl(\sup_{u\in[a,b]}\bigl|R_4^n(u)\bigr|
\Bigr)\leq K_{a,b}(k_n\Delta_n)^{1-\iota},
\end{eqnarray}
where $K_{a,b}$ depends only on $a$, and $b$ and is finite-valued.
\end{lemma}
\begin{pf} We start with
showing (\ref{prelim_b-1}). We define the set
%
\[
\label{proof_b_1} \mathcal{C}_i^n =\bigl \{ \bigl\llvert
\Delta_n^{-p/\beta}V_i^n(p)-\mu
_{p,\beta
}^{p/\beta}|\overline{\sigma}|_i^p
\bigr\rrvert > \tfrac{1}{2}\mu _{p,\beta
}^{p/\beta}|\overline{
\sigma}|_i^p \bigr\},\qquad i=k_n+3,\ldots ,n,
\]
and then we note that
%
\begin{eqnarray*}
\label{proof_b_2} %
1_{\{\mathcal{C}_i^n\}}&\leq&1 \bigl(\Delta_n^{-p/\beta}\bigl|V_i^n(p)
- \widehat{V}_i^n(p)\bigr|>\tfrac{1}{4}
\mu_{p,\beta}^{p/\beta}|\overline {\sigma }|_i^p
\bigr)
\\
&&{} +1 \bigl(\bigl|\Delta_n^{-p/\beta}\widehat {V}_i^n(p)-
\mu _{p,\beta}^{p/\beta}|\overline{\sigma}|_i^p\bigr|>
\tfrac{1}{4}\mu _{p,\beta
}^{p/\beta}|\overline{
\sigma}|_i^p \bigr).
\end{eqnarray*}
Hence we can apply (\ref{prelim_a-1}) and (\ref{prelim_a-2}) and conclude
%
\begin{equation}
\label{proof_b_3} \mathbb{E} \Bigl[\sup_{u\in\mathbb{R}_+}
\bigl(\bigl|r_i^1(u)\bigr|1_{\{\mathcal
{C}_i^n\}
}\bigr) \Bigr]\leq K
\bigl(\alpha_n\vee k_n^{- ({\beta
}/{(2p)}\wedge
{(\beta-p)}/{p} )+\iota} \bigr).
\end{equation}
We proceed with a sequence of inequalities. First, from Assumption \ref{assSB},
%
\begin{equation}
\label{proof_b_4} \mathbb{E}_{i-2}^n\biggl\llvert \int
_{(i-1)\Delta_n}^{i\Delta_n}(\alpha _u-\alpha
_{u-\Delta_n})\,du\biggr\rrvert \leq K\Delta_n^{1+{1}/{(r\vee1)}-\iota}.
\end{equation}
Next, if $\beta'<1$, we can decompose
%
\begin{equation}
\label{proof_b_5} \widehat{S}_t = \int_0^t
\int_{\mathbb{R}}x\mu_1(ds,dx)-t\int_{\mathbb
{R}}
\kappa(x)2\bigl|\nu'(x)\bigr|1_{\{\nu'(x)<0\}}\,dx,
\end{equation}
and separate accordingly $ \int_{(i-1)\Delta_n}^{i\Delta_n}\sigma
_{u-}\,d\widehat{S}_u$ and $ \int_{(i-2)\Delta_n}^{(i-1)\Delta
_n}\sigma
_{u-}\,d\widehat{S}_u$. For the difference of the integrals against time,
we can proceed exactly as in (\ref{proof_b_4}). Further, using the
algebraic inequality in (\ref{proof_a_3-a}), as well as Assumption~\ref{assA}
for the measure~$\nu'$, we have
%
\begin{equation}
\label{proof_b_6} \mathbb{E}_{i-1}^n\biggl\llvert \int
_{(i-1)\Delta_n}^{i\Delta_n}\int_{\mathbb
{R}}
\sigma_{u-}x\mu_1(du,dx) \biggr\rrvert ^{x}\leq
K\Delta_n^{x/\beta
'-\iota
}\qquad \mbox{for $x\leq\beta'$}.
\end{equation}
When $\beta'\geq1$, we can apply the Burkholder--Davis--Gundy
inequality and get
%
\begin{equation}
\label{proof_b_7} \mathbb{E}_{i-1}^n\biggl\llvert \int
_{(i-1)\Delta_n}^{i\Delta_n}\sigma _{u-}x\,d
\widehat{S}_u \biggr\rrvert ^{x}\leq K
\Delta_n^{x/\beta'-\iota
}\qquad \mbox{for $x\leq\beta'$}.
\end{equation}
The same inequalities hold for the analogous integrals involving
$\widetilde{S}$. Next, application of the Burkholder--Davis--Gundy and
H\"{o}lder inequalities, as well as Assumption~\ref{assSB} yields
%
\begin{equation}
\label{proof_b_8} \mathbb{E}_{i-2}^n\biggl\llvert \int
_{(i-1)\Delta_n}^{i\Delta_n}(\sigma _{u-}-
\sigma_{(i-2)\Delta_n-})\kappa(x)\widetilde{\mu}(du,dx) \biggr\rrvert \leq K
\Delta_n^{{2}/{\beta}-\iota}.
\end{equation}
Finally, denoting $\kappa'(x) = x-\kappa(x)$ and upon noting that
$\kappa'(x)$ is zero for $x$ sufficiently close to zero, we have
%
\begin{equation}\qquad
\label{proof_b_9} \mathbb{E}_{i-2}^n\biggl\llvert \int
_{(i-1)\Delta_n}^{i\Delta_n}(\sigma _{u-}-
\sigma_{(i-2)\Delta_n-})\kappa'(x)\mu(du,dx)\biggr\rrvert
^{\iota
}\leq K\Delta_n\qquad \forall\iota>0.
\end{equation}
Combining the estimates in (\ref{proof_b_4})--(\ref{proof_b_9}), as
well as the inequality $|\cos(x)-\cos(y)|\leq2|x-y|^p$ for every
$x,y\in\mathbb{R}$ and $p\in(0,1]$, we have
%
\begin{equation}\qquad
\label{proof_b_10} \mathbb{E} \Bigl[\sup_{u\geq a}
\bigl(\bigl|r_i^1(u)\bigr|1_{\{(\mathcal{C}_i^n)^c\}
}\bigr) \Bigr]\leq
K_a \bigl(\Delta_n^{{(\beta-\beta')}/{(\beta(\beta
'\vee
1))}-\iota}\vee
\Delta_n^{{1}/{\beta}\wedge{1}/{(r\vee
1)}-\iota} \bigr).
\end{equation}
Equations (\ref{proof_b_3}) and (\ref{proof_b_10}) yield (\ref{prelim_b-1}). We
continue next with (\ref{prelim_b-2}). This bound follows from a
first-order Taylor expansion of $f_{i,u}(x)$ and the bounds in (\ref
{decomp_1}) and~(\ref{prelim_a-1}).

We proceed by showing the result for $R_4^n(u)$. Using a second-order
Taylor expansion and the Cauchy--Schwarz inequality, as well as (\ref
{prelim_a-3-b}), we get
%
\begin{equation}
\label{proof_b_11} \mathbb{E} \Biggl(\sup_{u\in[a,b]}\Biggl\llvert
R_4^n(u) - \frac{\beta
}{p}e^{-C_{p,\beta}u^{\beta}}C_{p,\beta}u^{\beta}
\sum_{i=k_n+3}^n\widetilde{r}_i^4
\Biggr\rrvert \Biggr)\leq Kk_n,
\end{equation}
where
%
\[
\label{proof_b_12} \widetilde{r}_i^4 = \frac{|\sigma_{(i-2)\Delta_n-}|^p-|\overline
{\sigma
}|_i^p}{|\sigma_{(i-k_n-3)\Delta_n-}|^p}.
\]
Using (\ref{prelim_a-3}), we have
%
\begin{equation}
\label{proof_b_13} \mathbb{E}\Biggl\llvert \sum_{i=k_n+3}^n
\mathbb{E}_{i-k_n-3}^n\bigl(\widetilde {r}_i^4
\bigr)\Biggr\rrvert \leq K k_n.
\end{equation}
Further, without loss of generality (because $k_n\Delta_n\rightarrow
0$), we assume $n\geq2k_n+3$. Using the shorthand $\chi_i =
\widetilde
{r}_i^4-\mathbb{E}_{i-k_n-3}^n(\widetilde{r}_i^4)$, we then decompose
%
\begin{eqnarray*}
\label{proof_b_14} \sum_{i=k_n+3}^n
\chi_i &=& \sum_{j=1}^{k_n+1}A_j+
\sum_{i=2k_n+4+
(\lfloor{(n-k_n-2)}/{(k_n+1)}\rfloor-1 )(k_n+1)}^n\chi _i,
\\
\label{proof_b_15} A_j& =& \sum_{i=1}^{\lfloor{(n-k_n-2)}/{(k_n+1)}\rfloor}
\chi _{k_n+3+(j-1)+(i-1)(k_n+1)},\qquad j=1,\ldots,k_n+1.
\end{eqnarray*}
Applying the Burkholder--Davis--Gundy inequality for discrete
martingales and making use of (\ref{prelim_a-3-b}), we have
%
\begin{equation}
\label{proof_b_16} \mathbb{E}|A_j|\leq K(k_n
\Delta_n)^{-\iota},\qquad j=1,\ldots,k_n+1.
\end{equation}
Combining (\ref{proof_b_11}) and (\ref{proof_b_16}), we get the bound
in (\ref{prelim_b-4}).

We are left with (\ref{prelim_b-3}). The case $\beta/p\leq2$ follows from
%
\[
\mathbb{E}\bigl|r_i^3(u)\bigr|\leq K_{a,b}\bigl\llvert
\Delta_n^{-p/\beta}\widehat {V}_i^n(p)-
\mu_{p,\beta}^{p/\beta}\overline{V}_i^n(p)\bigr
\rrvert
\]
and by applying the bounds in (\ref{prelim_a-4-b})--(\ref
{prelim_a-4-c}). We now show (\ref{prelim_b-3}) for the case \mbox{$\beta/p>
2$}. We first decompose $r_i^3(u) = \sum_{j=1}^3\varrho_i^j(u)$, where
%
\begin{eqnarray*}
\label{proof_b_17} \varrho_i^1(u) &= & f'_{i,u}
\bigl(\Delta_n^{-p/\beta}\overline {V}_i^n(p)|
\overline{\sigma}|_i^p \bigr)\\
&&{}\times \Delta_n^{-p/\beta}
\bigl(\mu _{p,\beta}^{-p/\beta}\widehat{V}_i^n(p)
- |\overline{\sigma }|_i^p\overline{V}_i^n(p)
- \mu_{p,\beta}^{-p/\beta}\dot {V}_i^n(p)
\bigr),
\\
\label{proof_b_18} \varrho_i^2(u) &=& f'_{i,u}
(\widetilde{x} )\Delta _n^{-p/\beta
}\mu_{p,\beta}^{-p/\beta}
\dot{V}_i^n(p),
\\
\label{proof_b_19} %
\varrho_i^3(u) &=&
\bigl(f'_{i,u} (\widetilde{x} ) - f'_{i,u}
\bigl(\Delta_n^{-p/\beta}\overline{V}_i^n(p)|
\overline {\sigma }|_i^p \bigr) \bigr)
\\
&&{} \times\Delta_n^{-p/\beta} \bigl(\mu_{p,\beta}^{-p/\beta}
\widehat{V}_i^n(p) - |\overline {\sigma
}|_i^p\overline{V}_i^n(p) -
\mu_{p,\beta}^{-p/\beta}\dot {V}_i^n(p) \bigr)
\end{eqnarray*}
and $\widetilde{x}$ is a random number between $\Delta_n^{-p/\beta
}\mu
_{p,\beta}^{-p/\beta}\widehat{V}_i^n(p)$ and $\Delta_n^{-p/\beta
}|\overline{\sigma}|_i^p\overline{V}_i^n(p)$. We further introduce
%
\begin{eqnarray*}
\label{proof_b_20} \widetilde{\varrho}_i^1(u) =
\frac{G(p,u,\beta)}{|\sigma
_{(i-k_n-3)\Delta_n-}|^p} \Delta_n^{-p/\beta} \bigl(\mu_{p,\beta
}^{-p/\beta}
\widehat{V}_i^n(p) - |\overline{\sigma}|_i^p
\overline {V}_i^n(p) - \mu_{p,\beta}^{-p/\beta}
\dot{V}_i^n(p) \bigr)
\end{eqnarray*}
and note $G(p,u,\beta) = |\sigma_{(i-2)\Delta_n-}|^pf'_{i,u}
(|\sigma_{(i-2)\Delta_n-}|^p )$. Then direct calculation for the
function $xf'_{i,u}(x)$ and the boundedness of the process $|\sigma|$ yields
%
\[
\label{proof_b_21} %
\bigl|\varrho_i^1(u) - \widetilde{
\varrho}_i^1(u)\bigr|\leq K_{a,b}
\bigl(d_i^{(1)}+d_i^{(2)}
\bigr)e_i,
\]
where
%
\[
\label{proof_b_22} \cases{ %
d_i^{(1)} = \bigl|
\Delta_n^{-p/\beta}\overline {V}_i^n(p)-1\bigr|,
\vspace*{2pt}\cr
d_i^{(2)} = \bigl||\overline{\sigma}|_i^p
- |\sigma _{(i-2)\Delta_n-}|^p\bigr|+\bigl| |\sigma_{(i-2)\Delta_n-}|^p
- |\sigma _{(i-k_n-3)\Delta_n-}|^p\bigr|,
\vspace*{2pt}\cr
e_i = \Delta_n^{-p/\beta}\bigl\llvert \mu
_{p,\beta
}^{-p/\beta}\widehat{V}_i^n(p) - |
\overline{\sigma}|_i^p\overline {V}_i^n(p)
- \mu_{p,\beta}^{-p/\beta}\dot{V}_i^n(p)
\bigr\vert.}
\]
From here, we use the H\"{o}lder inequality and (\ref{prelim_a-2}),
(\ref{prelim_a-3-b}) and (\ref{prelim_a-4-b}) to get
%
\begin{equation}
\label{proof_b_23} \cases{ %
 \mathbb{E}\bigl|d_i^{(1)}e_i\bigr|
\leq K \bigl(\mathbb{E} \bigl[ \bigl(d_i^{(1)}
\bigr)^{{\beta
}/{(p+\beta\iota)}} \bigr] \bigr)^{{p}/{\beta}+\iota} \bigl(\mathbb {E}
\bigl(e_i^{{\beta}/{(\beta-p-\beta\iota)}} \bigr) \bigr)^{
{(\beta-p)}/{\beta}-\iota}
\vspace*{2pt}\cr
\hspace*{42pt}\leq Kk_n^{-1/2}(k_n\Delta_n)^{
{1}/{r}\wedge{(\beta-p)}/{\beta}-2\iota},
\vspace*{2pt}\cr
\mathbb {E}\bigl|d_i^{(2)}e_i\bigr|\leq\sqrt{
\mathbb{E}\bigl(d_i^{(2)}\bigr)^2\mathbb
{E}(e_i)^2}\leq K (k_n\Delta_n)^{1-\iota}.}
\end{equation}
For the sum $\sum_{i=k_n+3}^n\widetilde{\varrho}_i^1(u)$, using the
bounds in (\ref{prelim_a-4}) and (\ref{prelim_a-4-b}), we can proceed
exactly as for the analysis of $\sum_{i=k_n+3}^n\chi_i$ above and split
it into $k_n+1$ terms, which are the terminal values of discrete
martingales. Together, this yields
%
\begin{equation}
\label{proof_b_24} \mathbb{E} \Biggl(\sup_{u\in[a,b]}\Biggl\llvert \sum
_{i=k_n+3}^n\widetilde {\varrho
}_i^1(u)\Biggr\rrvert \Biggr)\leq K_{a,b}k_n(k_n
\Delta_n)^{-\iota}.
\end{equation}
Next, using the bound in (\ref{prelim_a-4-c}) as well as the
boundedness of the derivative $f_{i,u}'(x)$ (for $u\in[a,b]$), we have
%
\begin{equation}
\label{proof_b_25} \mathbb{E} \Bigl(\sup_{u\in[a,b]}\bigl|
\varrho_i^2(u) \bigr| \Bigr)\leq K_{a,b}
\Delta_n^{{(\beta-p)}/{\beta}\wedge{1}/{r}-\iota}.
\end{equation}
We continue with the term $\varrho_i^3(u)$. We first introduce the set
%
\[
\label{proof_b_26} \mathcal{E}_i^n = \bigl\{ \bigl\llvert
\mu_{p,\beta}^{-p/\beta}\widehat {V}_i^n(p) - |
\overline{\sigma}|_i^p\overline{V}_i^n(p)
- \mu _{p,\beta
}^{-p/\beta}\dot{V}_i^n(p)\bigr
\rrvert >1 \bigr\},\qquad i=k_n+3,\ldots,n.
\]
With this notation, using (\ref{prelim_a-4-b}) and the boundedness of
the derivative $f_{i,u}'(x)$ (for $u\in[a,b]$), we have
%
\begin{equation}
\label{proof_b_27} \mathbb{E} \Bigl(\sup_{u\in[a,b]}\bigl|
\varrho_i^3(u) \bigr|1_{\{\mathcal
{E}_i^n\}
} \Bigr)\leq
K_{a,b}(k_n\Delta_n)^{1-\iota}.
\end{equation}
Next using the boundedness of the second derivative $f_{i,u}^{\prime\prime}(x)$,
as well as the bounds in (\ref{prelim_a-4-b}) and (\ref{prelim_a-4-c}),
we get
%
\begin{equation}
\label{proof_b_28} \mathbb{E} \Bigl(\sup_{u\in[a,b]}\bigl|
\varrho_i^3(u)\bigr |1_{\{(\mathcal
{E}_i^n)^c\}} \Bigr)\leq
K_{a,b} \bigl((k_n\Delta_n)^{1-\iota}\vee
\Delta_n^{{(\beta-p)}/{\beta}\wedge{1}/{r}-\iota} \bigr).
\end{equation}
Combining (\ref{proof_b_23})--(\ref{proof_b_28}), we get the result in
(\ref{prelim_b-3}).
\end{pf}

\begin{lemma}\label{lema:prelim-c}
Under Assumptions \ref{assA} and \ref{assSB} and $k_n\asymp n^{\varpi}$ for $\varpi\in
(0,1)$, we have for $0<p<\beta$, $\iota>0$ arbitrarily small and every
$0<a<b<\infty$,
%
\begin{eqnarray}
\label{prelim_c-1} &&\frac{1}{n-k_n-2}\sup_{u\in[a,b]}\bigl|
\widehat{Z}_1^n(u) - \overline {Z}_1^n(u)
\bigr|
\nonumber
\\[-8pt]
\\[-8pt]
\nonumber
&&\qquad= o_p \bigl(\alpha_n\vee k_n^{- ({\beta
}/{(2p)}\wedge
{(\beta-p)}/{p} )+\iota}
\vee\sqrt{\Delta_n} \bigr),
\end{eqnarray}
and further if $p<\beta/2$,
%
\begin{eqnarray}
\label{prelim_c-2} %
&&\frac{1}{n-k_n-2}\sup_{u\in[a,b]}\bigl|
\widehat{Z}_2^n(u) - \overline {Z}_2^{(a,n)}(u)
- \overline{Z}_2^{(b,n)}(u)\bigr|
\nonumber
\\[-8pt]
\\[-8pt]
\nonumber
&&\qquad = o_p \bigl( k_n^{-({1}/{2}) ({\beta}/{p}\wedge3
)+\iota
}\vee
k_n^{-1/2}(k_n\Delta_n)^{{1}/{r}\wedge{(\beta
-p)}/{\beta
}-\iota}
\bigr). %
\end{eqnarray}
\end{lemma}
%
\begin{pf}
We start with
(\ref{prelim_c-1}). We split $\widehat{Z}_1^n(u) - \overline{Z}_1^n(u)
= E_1^n(u)+E_2^n(u)$ with $E_1^n(u) = \sum_{i=k_n+3}^n(z_{i}^1(u) -
\overline{z}_{i}^1(u))1_{\{\mathcal{C}_i^n\}}$ and $E_2^n(u) = \sum_{i=k_n+3}^n(z_{i}^1(u) - \overline{z}_{i}^1(u))1_{\{(\mathcal
{C}_i^n)^c\}}$. For $E_1^n(u)$, using Lemma~\ref{lema:prelim-a}, we
easily have
%
\begin{equation}
\label{proof_c_2} \frac{1}{n-k_n-2}\mathbb{E} \Bigl(\sup_{u\in[a,b]}\bigl|E_1^n(u)
\bigr| \Bigr)\leq K_{a,b} \bigl(\alpha_n\vee
k_n^{- ({\beta}/{(2p)}\wedge{(\beta-p)}/{p} )+\iota} \bigr).
\end{equation}

We proceed with $E_2^n(u)$. We first note that
%
\begin{equation}
\label{proof_c_3} \mathbb{E}_{i-2}^n \bigl[
\bigl(z_{i}^1(u) - \overline{z}_{i}^1(u)
\bigr)1_{\{
(\mathcal{C}_i^n)^c\}} \bigr] = 0.
\end{equation}
Further, using the algebraic inequalities $|\cos(x)-\cos(y)|^2\leq
2|x-y|$ for $x,y\in\mathbb{R}$ and $|e^{-x}-e^{-y}|^2\leq2|x-y|$ for
$x,y\in\mathbb{R}_+$, as well as the definition of the set $\mathcal
{C}_i^n$, we get
%
\[
\mathbb{E}_{i-2}^n\bigl\llvert \bigl(z_{i}^1(u)
- \overline{z}_{i}^1(u)\bigr)1_{\{
(\mathcal{C}_i^n)^c\}} \bigr
\rrvert ^2\leq K_{a,b}\bigl\llvert \Delta
_n^{-p/\beta
}V_i^n(p) -
\mu_{p,\beta}^{p/\beta}|\sigma_{(i-2)\Delta
_n-}|^p\bigr
\rrvert .
\]
Applying the above two inequalities, the bounds in (\ref{prelim_a-1}),
(\ref{prelim_a-2}) and (\ref{prelim_a-3-b}), as well as the algebraic
inequality $2xy\leq x^2+y^2$ for $x,y\in\mathbb{R}$, we have
%
\begin{eqnarray*}
\label{proof_c_4} %
\mathbb{E} \bigl( E_2^n(u)
\bigr)^2 &=& \mathbb{E} \Biggl(\sum_{i=k_n+3}^n
\bigl(z_{i}^1(u) - \overline{z}_{i}^1(u)
\bigr)^21_{\{(\mathcal
{C}_i^n)^c\}} \Biggr)
\\
&&{} +\mathbb{E} \biggl(\sum_{i,j\dvtx
|i-j|=1}\bigl(z_{i}^1(u)
- \overline{z}_{i}^1(u)\bigr)1_{\{(\mathcal{C}_i^n)^c\}
}
\bigl(z_{j}^1(u) - \overline{z}_{j}^1(u)
\bigr)1_{\{(\mathcal{C}_j^n)^c\}
} \biggr)
\\
& \leq& K_{a,b}\sum_{i=k_n+3}^n
\mathbb{E}\bigl\llvert \Delta_n^{-p/\beta
}V_i^n(p)
- \mu_{p,\beta}^{p/\beta}|\sigma_{(i-2)\Delta
_n-}|^p\bigr
\rrvert
\\
&\leq& K_{a,b}\Delta_n^{-1} \bigl(
\alpha_n\vee k_n^{- ({\beta}/{(2p)}\wedge{(\beta-p)}/{p} )+\iota}\vee(k_n\Delta
_n)^{{1}/{r}\wedge1-\iota} \bigr).
\end{eqnarray*}
As a result, $\frac{1}{\sqrt{n-k_n-2}}E_2^n(u) \stackrel{\mathbb
{P}}{\longrightarrow} 0$ finite-dimensionally in $u$. Finally, we need
to show that the convergence holds uniformly in $u\in[a,b]$. For this
we apply a criteria for tightness on the space of continuous functions
equipped with the uniform topology; see, for example, Theorem~12.3 of
\cite{Billingsley}. Using again (\ref{proof_c_3}), we have
%
\begin{eqnarray*}
\label{proof_c_5} %
&&\mathbb{E} \bigl( E_2^n(u)-E_2^n(v)
\bigr)^2
\\
&&\qquad \leq K\mathbb {E} \Biggl(\sum_{i=k_n+3}^n
\bigl(z_{i}^1(u) - \overline{z}_{i}^1(u)
- z_{i}^1(v) + \overline{z}_{i}^1(v)
\bigr)^21_{\{(\mathcal{C}_i^n)^c\}} \Biggr).
\end{eqnarray*}
Hence for arbitrarily small $\iota>0$,
%
\begin{eqnarray*}
\label{proof_c_6} %
&&\frac{1}{n-k_n-2}\mathbb{E} \Biggl(\sum
_{i=k_n+3}^n\bigl(z_{i}^1(u) -
\overline{z}_{i}^1(u) - z_{i}^1(v)
+ \overline{z}_{i}^1(v)\bigr)^21_{\{
(\mathcal{C}_i^n)^c\}}
\Biggr)
\\
&&\qquad \leq K \bigl\{\bigl|u^{\beta
}-v^{\beta}\bigr|^2
\vee|u-v|^{\beta-\iota} \bigr\},
\end{eqnarray*}
and since $\beta>1$, we have $\frac{1}{\sqrt{n-k_n-2}}\sup_{u\in
[a,b]}|E_2^n(u)| \stackrel{\mathbb{P}}{\longrightarrow} 0$. We turn
next to~(\ref{prelim_c-2}). We first introduce some additional
notation. Based on a second-order Taylor expansion of the function
$f_{i,u}(x)$, we can further decompose $\widehat{Z}_2^n(u) = \widehat
{Z}_2^{(a,n)}(u) +\widehat{Z}_2^{(b,n)}(u)+\widehat{Z}_2^{(c,n)}(u)$,
with $\widehat{Z}_2^{(k,n)}(u) = \sum_{i=k_n+3}^nz_i^{(k,2)}(u)$ for
$k=a,b,c$, where $z_i^{(c,2)}(u) =
z_i^{2}(u)-z_i^{(a,2)}(u)-z_i^{(b,2)}(u)$ and
%
\begin{eqnarray*}
\label{proof_c_7} z_i^{(a,2)}(u)& =& f_{i,u}'
\bigl( |\overline{\sigma}|_i^p \bigr)|\overline{
\sigma}|_i^p \bigl( \Delta_n^{-p/\beta}
\overline {V}_i^n(p) - 1 \bigr),
\\
\label{proof_c_8} z_i^{(b,2)}(u) &=& \tfrac{1}{2}f_{i,u}^{\prime\prime}
\bigl( |\overline{\sigma}|_i^p \bigr) \bigl(|\overline{
\sigma}|_i^p\bigr)^2 \bigl(
\Delta_n^{-p/\beta
}\overline {V}_i^n(p) -
1 \bigr)^2.
\end{eqnarray*}
Note further that
%
\[
\label{proof_c_9} \cases{ %
|\sigma_{(i-2)\Delta_n-}|^pf_{i,u}'
\bigl( |\sigma _{(i-2)\Delta_n-}|^p \bigr) = G(p,u,\beta),
\vspace*{2pt}\cr
|\sigma _{(i-2)\Delta
_n-}|^{2p}f_{i,u}^{\prime\prime} \bigl(
|\sigma_{(i-2)\Delta_n-}|^p \bigr) = H(p,u,\beta). }
\]
Direct calculation, and using the boundedness of the process $\sigma$
by Assumption~\ref{assSB}, shows
%
\begin{eqnarray*}
\label{proof_c_10} %
&& \bigl\llvert |\overline{\sigma}|_i^pf_{i,u}'
\bigl( |\overline{\sigma}|_i^p \bigr) - G(p,u,\beta)\bigr
\rrvert + \bigl\llvert \bigl(|\overline{\sigma}|_i^{p}
\bigr)^2f_{i,u}^{\prime\prime} \bigl( |\overline {\sigma
}|_i^p \bigr) - H(p,u,\beta)\bigr\rrvert
\\
& &\qquad\leq K_{a,b}\bigl\llvert |\overline{\sigma}|_i^p-|
\sigma_{(i-2)\Delta_n-}|^p\bigr\rrvert ,\qquad u\in [a,b],
i=k_n+3,\ldots,n,
\end{eqnarray*}
for some finite-valued constant $K_{a,b}$ which depends only $a$ and
$b$. From here, using the bounds in (\ref{prelim_a-2}) and (\ref
{prelim_a-3-b}), we have
%
\[
\label{proof_c_11} \mathbb{E} \Bigl(\sup_{u\in[a,b]}\bigl|z_i^{(a,2)}(u)
- \overline {z}_i^{(a,2)}(u)\bigr| \Bigr)\leq K_{a,b}
\bigl( k_n^{-1/2}(k_n\Delta _n)^{{1}/{r}\wedge{(\beta-p)}/{\beta}-\iota}
\bigr),
\]
and similarly
%
\[
\label{proof_c_11_b} \mathbb{E} \Bigl(\sup_{u\in[a,b]}\bigl|z_i^{(b,2)}(u)
- \overline {z}_i^{(b,2)}(u)\bigr| \Bigr)\leq K_{a,b}
\bigl(k_n^{-{\beta
}/{(2p)}+\iota} \vee k_n^{-1/2}(k_n
\Delta_n)^{{1}/{r}\wedge{(\beta-p)}/{\beta}-\iota} \bigr).
\]
Therefore,
%
\begin{eqnarray}
\label{proof_c_12} %
&&\frac{1}{n-k_n-2}\sup_{u\in[a,b]}\bigl
\llvert \widehat{Z}_2^{(a,n)}(u) - \overline{Z}_2^{(a,n)}(u)
+ \widehat{Z}_2^{(b,n)}(u) - \overline {Z}_2^{(b,n)}(u)
\bigr\rrvert
\nonumber
\\[-8pt]
\\[-8pt]
\nonumber
& &\qquad= o_p \bigl( k_n^{-{\beta
}/{(2p)}+\iota} \vee
k_n^{-1/2}(k_n\Delta_n)^{{1}/{r}\wedge
{(\beta
-p)}/{\beta}-\iota}
\bigr). %
\end{eqnarray}
We are left with $\widehat{Z}_2^{(c,n)}(u)$. Using the boundedness of
the derivatives in (\ref{decomp_1}), we have
%
\[
\label{proof_c_13} \bigl\llvert z_i^{(c,2)}(u)\bigr\rrvert \leq
K_{a,b}\bigl\llvert \Delta_n^{-p/\beta
}
\overline{V}_i^n(p) - 1\bigr\rrvert ^{x},\qquad 2<x<
\beta/p\wedge3.
\]
From here, applying (\ref{prelim_b-2}), we have
%
\begin{equation}
\label{proof_c_14} \frac{1}{n-k_n-2}\sup_{u\in[a,b]}\bigl\llvert
\widehat {Z}_2^{(c,n)}(u)\bigr\rrvert = o_p
\bigl(k_n^{-({1}/{2}) ({\beta}/{p}\wedge3
)+\iota
} \bigr).
\end{equation}
Combining the results in (\ref{proof_c_12}) and (\ref{proof_c_14}), we
get (\ref{prelim_c-2}).
\end{pf}

\begin{lemma}\label{lema:prelim-d}
Let $p\in(0,\beta/2)$. If $k_n\asymp n^{\varpi}$ for $\varpi\in(0,1)$,
we have
%
\begin{equation}
\label{prelim_d-1} \frac{1}{\sqrt{n-k_n-2}}\pmatrix{
\overline{Z}_1^n(\mathbf{u})
\vspace*{2pt}\cr
\overline {Z}_2^{(a,n)}(\mathbf{u})} \stackrel{\mathcal{L}} {\longrightarrow} \zeta(\mathbf{u}),
\end{equation}
where $\zeta(\mathbf{u})$ is a Gaussian process with covariance
function given by
%
\begin{equation}
\label{prelim_d-2}\pmatrix{ 1
\vspace*{2pt}\cr
G(p,u,\beta) }'\overline{
\Xi}(p,u,v,\beta)\pmatrix{ 1
\vspace*{2pt}\cr
G(p,v,\beta) }, \qquad u,v\in\mathbb{R}_+,
\end{equation}
for $\overline{\Xi}(p,u,v,\beta) = \Xi_0(p,u,v,\beta)+2\Xi
_1(p,u,v,\beta
)$. The convergence in (\ref{prelim_d-1}) is in the space of continuous
functions $\mathbb{R}_+\rightarrow\mathbb{R}^2$ equipped with the local
uniform topology. The convergence result for $\overline{Z}_1^n(\mathbf
{u})$ in (\ref{prelim_d-1}) continues to hold for $p\in[\beta
/2,\beta)$.

Further, for some $\iota>0$,
%
\begin{eqnarray}
\label{prelim_d-3} %
&&\frac{k_n}{n-k_n-2}\overline{Z}_2^{(b,n)}(u)\nonumber\\
&&\quad{}-
\frac
{1}{2}H(p,u,\beta ) \bigl( \Xi_0^{(2,2)}(p,u,u,
\beta) +2\Xi_1^{(2,2)}(p,u,u,\beta) \bigr)
\\
&&\qquad = o_p \bigl((k_n\Delta_n)^{1-{2p}/{\beta}\vee{1}/{2}-\iota}
\bigr),\nonumber %
\end{eqnarray}
locally uniformly in $u\in\mathbb{R}_+$.
\end{lemma}
%
\begin{pf}
We can write
%
\begin{equation}
\label{proof_d_1}\pmatrix{\overline{Z}_1^n(u)
\vspace*{2pt}\cr
\overline{Z}_2^{(a,n)}(u) } = \sum_{i=k_n+1}^{n-k_n-1}\bolds{
\zeta}_i(u)+E_l(u)+E_r(u),
\end{equation}
where
%
\begin{eqnarray*}
\label{proof_d_2} \bolds{\zeta}_i(u) &=& \pmatrix{\cos \bigl( u\Delta_n^{-1/\beta}
\mu_{p,\beta
}^{-1/\beta
} \bigl(\Delta_i^nS-
\Delta_{i-1}^nS\bigr) \bigr) - \mathcal{L}(p,u,\beta)
\vspace*{2pt}\cr
G(p,u,\beta) \bigl[\Delta_n^{-p/\beta} \mu_{p,\beta}^{-p/\beta
}\bigl|
\Delta _i^nS-\Delta_{i-1}^nS\bigr|^p-1
\bigr] },
\\
\label{proof_d_3} E_l(u)& =& \sum_{i=2}^{k_n}
\frac{i-1}{k_n}\pmatrix{ 0
\vspace*{2pt}\cr
\bolds{\zeta}^{(2)}_i(u) }-\sum_{i=k_n+1}^{k_n+2}\pmatrix{ \bolds{\zeta}^{(1)}_i(u)
\vspace*{2pt}\cr
0 },
\\
\label{proof_d_4} E_r(u)& =& \sum_{i=n-k_n}^{n-2}
\frac{n-1-i}{k_n}\pmatrix{ 0
\vspace*{2pt}\cr
\bolds{\zeta}^{(2)}_i(u)}+\sum_{i=n-k_n}^{n}\pmatrix{ \bolds{\zeta}^{(1)}_i(u)
\vspace*{2pt}\cr
0 }.
\end{eqnarray*}
We note that for $u\in\mathbb{R}_+$,
%
\begin{equation}
\label{proof_d_5} \mathbb{E}_{i-2}^n \bigl(\bolds{
\zeta}_i(u) \bigr) = 0,\qquad i=2,\ldots,n.
\end{equation}
Further, making using of the inequality $|\cos(x)-\cos(y)|\leq
2|x-y|^p$ for every $p\in(0,1]$ and $x,y\in\mathbb{R}$, we have for
$u,v\in\mathbb{R}_+$,
%
\begin{equation}\qquad
\label{proof_d_6} \mathbb{E}_{i-2}^n \bigl(\bolds{
\zeta}^{(1)}_i(u)-\bolds{\zeta }^{(1)}_i(v)
\bigr)^2\leq K|u-v|^p\vee\bigl|u^{\beta}-v^{\beta}\bigr|^2,\qquad
1<p<\beta.
\end{equation}
Making use of (\ref{proof_d_5}) and the fact that $\bolds{\zeta
}^{(2)}_i(u)$ depends on $u$ only through $H(p,u,\beta)$ and $\sup_{u\in
\mathbb{R}^+}|H(p,u,\beta)|$ is a finite constant, we have
%
\begin{equation}
\label{proof_d_7} \frac{1}{k_n}\mathbb{E} \Bigl(\sup_{u\in\mathbb
{R}_+}\bigl|E_l(u)\bigr|^2
\Bigr)\leq K.
\end{equation}
Making use of (\ref{proof_d_6}) and the differentiability of
$G(p,u,\beta)$ in $u$, we also have
%
\[
\label{proof_d_8} \frac{1}{k_n}\mathbb{E} \bigl(E_r(u)-E_r(v)
\bigr)^2\leq \bigl|F(u)-F(v)\bigr|^p,
\]
for some increasing function $F(\cdot)$ and some $p>1$. Applying then a
criteria for tightness on the space of continuous functions equipped
with the uniform topology (see, e.g., Theorem~12.3 in \cite
{Billingsley}) as well as making use of the fact that $k_n\Delta
_n\rightarrow0$, we have locally uniformly in $u$,
%
\begin{equation}
\label{proof_d_9} \frac{1}{\sqrt{n-k_n-2}}E_r(u) \stackrel{\mathbb {P}} {
\longrightarrow} 0.
\end{equation}
We are left with the first term on the right-hand side of (\ref
{proof_d_1}). First, we establish convergence for this term
finite-dimensionally in $u$. We have the decomposition
%
\[
\label{proof_d_10} \sum_{i=k_n+1}^{n-k_n-1}\bolds{
\zeta}_i(u) = \sum_{i=k_n+1}^{n-k_n-1}
\bigl(\bolds{\zeta}_i(u)-\mathbb{E}_{i-1}^n
\bigl(\bolds {\zeta}_i(u)\bigr)\bigr)+\sum
_{i=k_n}^{n-k_n-2}\mathbb{E}_i^n
\bigl(\bolds{\zeta }_{i+1}(u)\bigr).
\]
From here, we can apply a c.l.t. for triangular arrays (see, e.g., Theorem~2.2.13 of \cite{JP}) to establish that $\frac{1}{\sqrt
{n-2k_n-1}}\sum_{i=k_n+1}^{n-k_n-1}\bolds{\zeta}_i(u)$
converges finite-dimensionally in $u$ to $\zeta(\mathbf{u})$. This
convergence holds also locally uniformly in $u$ using the bound in
(\ref{proof_d_6}) and Theorem VI.4.1 in \cite{JS}. Combining the latter
with the asymptotic negligibility results in (\ref{proof_d_7}) and
(\ref{proof_d_9}), together with the fact that $k_n/n\rightarrow0$, we
have the result in (\ref{prelim_d-1}). Furthermore, since $\overline
{Z}_1^n(\mathbf{u})$ depends on $p$ only through $\mu_{p,\beta}$, the
marginal convergence in (\ref{prelim_d-1}) involving $\overline
{Z}_1^n(\mathbf{u})$ holds for any $p\in(0,\beta)$.

We turn next to (\ref{prelim_d-3}). We denote
%
\[
\label{proof_d_11} \chi_i = k_n \bigl(\Delta_n^{-p/\beta}
\overline{V}_i^n(p) -1 \bigr)^2 - \bigl(
\Xi_0^{(2,2)}(p,u,u,\beta) +2\Xi_1^{(2,2)}(p,u,u,
\beta) \bigr),
\]
and we note that $\Xi_0^{(2,2)}(p,u,u,\beta) $ and $\Xi
_1^{(2,2)}(p,u,u,\beta) $ do not depend on $u$.

Without loss of generality we can assume $n\geq2k_n+3$, and then we set
%
\[
\label{proof_d_12} A_j = \sum_{i=1}^{\lfloor{(n-k_n-2)}/{(k_n+1)}\rfloor}
\chi _{k_n+3+(j-1)+(i-1)(k_n+1)},\qquad j=1,\ldots,k_n+1.
\]
Since $\mathbb{E}|\chi_i|<K$,
%
\begin{equation}
\label{proof_d_13} \Biggl\llvert \sum_{i=k_n+3}^n
\chi_i - \sum_{j=1}^{k_n+1}A_j
\Biggr\rrvert = O_p(k_n).
\end{equation}
Further, direct computation shows
%
\[
\label{proof_d_14} \mathbb{E}_{i-k_n-3}^n (\chi_i )=
0,\qquad i=k_n+3,\ldots,n,
\]
and applying the Burkholder--Davis--Gundy inequality for discrete
martingales, we have
%
\begin{equation}
\label{proof_d_15} \mathbb{E}|A_j|^x\leq K(k_n
\Delta_n)^{- ({x}/{2}\vee
1
)}, \qquad 1\leq x<\frac{\beta}{2p}.
\end{equation}
Using inequality in means we further have
%
\[
\label{proof_d_16} \Biggl\llvert \frac{1}{k_n+1}\sum
_{j=1}^{k_n+1}A_j\Biggr\rrvert
^x\leq\frac
{1}{k_n+1}\sum_{j=1}^{k_n+1}|A_j|^x,\qquad
1\leq x<\frac{\beta
}{2p}.
\]
Applying the above inequality with $x$ sufficiently close to $\beta
/(2p)$ and the bound in (\ref{proof_d_15}), we have $ \Delta
_n(k_n\Delta
_n)^{{2p}/{\beta}\wedge{1}/{2}-1+\iota}
\sum_{j=1}^{k_n+1}A_j \stackrel{\mathbb{P}}{\longrightarrow} 0$, and
together with the result in (\ref{proof_d_13}), this implies (\ref
{prelim_d-3}).
\end{pf}

\subsection{Proofs of Theorems \texorpdfstring{\protect\ref{thm:cp}}{1} and 
\texorpdfstring{\protect\ref{thm:clt}}{2}}
Theorem~\ref{thm:cp} and (\ref{clt_3}) of Theorem~\ref{thm:clt} follow
readily by combining Lemmas \ref{lema:prelim-a}--\ref{lema:prelim-d}
[and using (\ref{prelim_a-2}) for bounding $\widehat{Z}_2^n(u)$ in the
proof of Theorem~\ref{thm:cp}]. To show (\ref{clt_5}), we note first
that $H(p,u,\beta)$ and $\Xi_i(p,u,u,\beta)$, for $i=0,1$, are
continuously differentiable in $\beta$. For $H(p,u,\beta)$ this is
directly verifiable, and for $\Xi_i(p,u,u,\beta)$ with $i=0,1$, this
follows from the continuous differentiability of the characteristic
function $\beta\rightarrow e^{-A_{\beta}u^{\beta}}$ for $u\in
\mathbb
{R}_+$. Moreover, the derivative $\nabla_{\beta}H(p,u,\beta)$ is
bounded in $u$. From here, (\ref{clt_5}) follows from an application of
the continuous mapping theorem.

\subsection{Proof of Theorem \texorpdfstring{\protect\ref{thm:gmm}}{3}}

We denote the true value of the parameter $\beta$ with $\beta_0$. Then
the claim in (\ref{gmm_6}) will follow if we can show the following:
%
\begin{equation}
\label{proof_gmm_1} \widehat{\mathbf{m}}\bigl(p,\widehat{\mathbf{u}},\widehat{\beta
}^{fs},\mathbf {u},\beta\bigr) \stackrel{\mathbb{P}} {\longrightarrow}
\mathbf {m}(p,\mathbf {u},\beta) \qquad\mbox{uniformly in $\beta\in[1,2]$,}
\end{equation}
where $\mathbf{m}(p,\mathbf{u},\beta)$ is defined via
%
\begin{eqnarray}
\mathbf{m}(p,\mathbf{u},\beta)_i &=& \int_{u_l^i}^{u_h^i}
\bigl(\log \bigl(\mathcal{L}(p,u,\beta_0)\bigr) - \log\bigl(\mathcal{L}(p,u,
\beta)\bigr) \bigr)\,du,
\nonumber
\\
\qquad\quad\label{proof_gmm_2} \sqrt{n} \widehat{\mathbf{m}}\bigl(p,\widehat{\mathbf{u}},
\widehat {\beta }^{fs},\mathbf{u},\beta_0\bigr) &\stackrel{
\mathcal{L}} {\longrightarrow }& \mathbf{W}^{1/2}(p,\mathbf{u},
\beta_0)\times\mathcal{\mathbf{N}},
\end{eqnarray}
where $\mathcal{\mathbf{N}}$ is $K\times1$ standard normal vector and
%
\begin{equation}\qquad
\label{proof_gmm_3} \mathbf{M}(p,\widehat{\mathbf{u}},\beta) \stackrel{\mathbb {P}} {
\longrightarrow} \mathbf{M}(p,\mathbf{u},\beta) \qquad\mbox {uniformly in a
neighborhood of $\beta_0$.}
\end{equation}
This is because $\mathbf{m}(p,\mathbf{u},\beta) = \mathbf{0}$ if and
only if $\beta=\beta_0$ and $W(p,\mathbf{u},\beta_0)$ is positive definite.

 We start with (\ref{proof_gmm_1}). We have
 \[
 \int_{\widehat
{u}_l^i}^{\widehat{u}_h^i} \log\bigl(\mathcal{L}(p,u,\beta)\bigr)\,du \stackrel
{\mathbb{P}}{\longrightarrow} \int_{u_l^i}^{u_h^i} \log\bigl(\mathcal
{L}(p,u,\beta)\bigr)\,du
\]
uniformly in $\beta\in[1,2]$ for $i=1,\ldots,K$ because
of $\widehat{\mathbf{u}}_l \stackrel{\mathbb{P}}{\longrightarrow
} \mathbf{u}_l$ and $\widehat{\mathbf{u}}_h \stackrel{\mathbb
{P}}{\longrightarrow} \mathbf{u}_h$ as well as the continuity of the
function $u^{\beta}$ in $\beta$ for every $u\in\mathbb{R}_+$, and the
argument can be used to show (\ref{proof_gmm_3}). To show (\ref
{proof_gmm_1}) it remains to show $\int_{\widehat{u}_l^i}^{\widehat
{u}_h^i} \log(\widehat{\mathcal{L}}^n(p,u,\widehat{\beta
}^{fs})')\,du \stackrel{\mathbb{P}}{\longrightarrow} \int_{u_l^i}^{u_h^i}
\log(\mathcal{L}(p,u,\beta_0))\,du$ for $i=1,\ldots,K$.\vspace*{-3pt} Due the continuous
differentiability of the de-biasing\vspace*{1pt} term in $\beta$, $\widehat{\beta
}^{fs} \stackrel{\mathbb{P}}{\longrightarrow} \beta_0$ and the
asymptotic boundedness of $\widehat{\mathbf{u}}_l$ and $\widehat
{\mathbf
{u}}_h$ and of $\widehat{\mathcal{L}}^n(p,u,\widehat{\beta}^{fs})'$
from below, we have $\int_{\widehat{u}_l^i}^{\widehat{u}_h^i} [\log
(\widehat{\mathcal{L}}^n(p,u,\widehat{\beta}^{fs})')- \log
(\widehat
{\mathcal{L}}^n(p,u,\beta_0))]\,du \stackrel{\mathbb
{P}}{\longrightarrow
} 0$. From here (\ref{proof_gmm_1}) follows by applying Theorem~\ref{thm:cp}.

We are left with (\ref{proof_gmm_2}). This result follows from applying
the uniform convergence of $\widehat{\mathcal{L}}^n(p,u,\widehat
{\beta
}^{fs})'$ in Theorem~\ref{thm:clt}.

Finally, (\ref{gmm_7}) follows from the continuity of $G(p,\mathbf
{u},\beta)$ and $W^{-1}(p,\mathbf{u},\beta)$ in $\mathbf{u}$ and
$\beta$.

\subsection{Proof of Theorem \texorpdfstring{\protect\ref{thm:cont}}{4}}
We will use the shorthand notation $v_n = \rho u_n$. We start with the
following lemma.

\begin{lemma}\label{lema:prelim-e}
Under the conditions of Theorem~\ref{thm:cont} we have
%
\begin{eqnarray}
\label{prelim_e-1} \widehat{\mathcal{L}}^n\bigl(p,u_n,
\widehat{\beta}^{fs}\bigr)' - \mathcal {L}
\bigl(p,u_n,\widehat{\beta}^{fs}\bigr) &=& O_p
\bigl(\sqrt{\Delta_n}u_n^{2}\bigr),
\\
\label{prelim_e-2} %
\frac{\sqrt{n}}{u_n^2-v_n^2}\widehat{Z}_n &\stackrel{
\mathcal {L}} {\longrightarrow} &\frac{1}{24C_{p,2}}Z_1-
\frac{2}{p}C_{p,2}Z_2, %
\end{eqnarray}
where
\begin{eqnarray*}
\widehat{Z}_n &=& \frac{1}{C_{p,2}u_n^2} \bigl(\widehat{\mathcal
{L}}^n\bigl(p,u_n,\widehat{\beta}^{fs}
\bigr)' - \mathcal{L}\bigl(p,u_n,\widehat {\beta
}^{fs}\bigr) \bigr) \\
&&{}- \frac{1}{C_{p,2}v_n^2} \bigl(\widehat{\mathcal
{L}}^n\bigl(p,v_n,\widehat{\beta}^{fs}
\bigr)' - \mathcal{L}\bigl(p,v_n,\widehat {\beta
}^{fs}\bigr) \bigr).
\end{eqnarray*}
\end{lemma}

\begin{pf}
We use the same
decomposition of $\widehat{\mathcal{L}}^n(p,u,\beta) - \mathcal
{L}(p,u,\beta)$ as in the proofs of Theorems \ref{thm:cp} and \ref
{thm:clt}. We start with the leading terms $\overline{Z}_1^n(u_n)$,
$\overline{Z}_2^{(a,n)}(u_n)$ and $\overline{Z}_2^{(b,n)}(u_n)$. Using
Taylor's series expansion, we have for any $u\in\mathbb{R}_+$ and $Z\in
\mathbb{R}$,
%
\begin{eqnarray*}
\label{proof_e_1} \cos(uZ)-1 &=& -\frac{u^2Z^2}{2}+\frac{u^4Z^4}{24}+R(uZ),\qquad
\bigl|R(uZ)\bigr|\leq K|uZ|^6,
\\
\label{proof_e_2} 1-e^{-u^2} &=& u^2-\frac{u^4}{2}+O
\bigl(u^6\bigr) \qquad\mbox{as $u\rightarrow 0$}.
\end{eqnarray*}
Using this approximation we have (note that when $L_t$ is a Brownian
motion, then $A_{\beta} = 1$ and so $C_{p,\beta} = 1/\mu_{p,\beta}$)
%
\begin{eqnarray}
\label{proof_e_3} %
&&\frac{1}{C_{p,2}u_n^2}\overline{Z}_1^n(u_n)
- \frac
{1}{C_{p,2}v_n^2}\overline{Z}_1^n(v_n)
\nonumber
\\[-8pt]
\\[-8pt]
\nonumber
&&\qquad = \frac
{u_n^2-v_n^2}{24C_{p,2}}\sum_{i=k_n+3}^n
\biggl[\frac{n^2(\Delta
_i^nS-\Delta_{i-1}^nS)^4}{\mu_{p,2}^{2}}-\frac{12}{\mu
_{p,2}^{2}} \biggr]+O_p
\bigl(u_n^4\sqrt{n}\bigr). %
\end{eqnarray}
We similarly get
%
\begin{eqnarray}
\label{proof_e_4} %
&&\frac{1}{C_{p,2}u_n^2}\overline{Z}_2^{(a,n)}(u_n)
- \frac
{1}{C_{p,2}v_n^2}\overline{Z}_2^{(a,n)}(v_n)
\nonumber
\\[-8pt]
\\[-8pt]
\nonumber
&&\qquad = \bigl(v_n^2-u_n^2\bigr)
\frac{2}{p}C_{p,2}\sum_{i=k_n+3}^n
\bigl(\Delta _n^{-p/2}\overline{V}_i^n(p)-1
\bigr)+O_p\bigl(u_n^4\sqrt{n}\bigr),
\end{eqnarray}
and also
%
\begin{eqnarray}
\label{proof_e_5} %
&&\frac{k_n}{n-k_n-2}\sum_{i=k_n+3}^n
\bigl(\Delta_n^{-p/2}\overline {V}_i^n(p)-1
\bigr)^2 \nonumber\\
&&\quad{}- \bigl(\Xi_0^{(2,2)}(p,u_n,u_n,2)+2
\Xi _1^{(2,2)}(p,u_n,u_n,2) \bigr)
\\
&& \qquad= O_p (\sqrt{k_n\Delta _n} ). \nonumber%
\end{eqnarray}
As in Lemma~\ref{lema:prelim-d}, it is easy to show
%
\begin{equation}
\label{proof_e_6} \frac{1}{\sqrt{n-k_n-2}}\sum_{i=k_n+3}^n
\pmatrix{\displaystyle\frac{n^2(\Delta_i^nS-\Delta_{i-1}^nS)^4}{\mu
_{p,2}^{2}}-\frac{12}{\mu_{p,2}^{2}}
\vspace*{2pt}\cr
\displaystyle\Delta_n^{-p/2}\overline {V}_i^n(p)-1
} \stackrel{\mathcal{L}} {\longrightarrow} \pmatrix{ Z_1
\vspace*{2pt}\cr
Z_2 }.
\end{equation}
Next, using Taylor's expansion as well as $\widehat{\beta}^{fs}-2 =
o_p(k_nu_n^{2}\sqrt{\Delta_n})$, we have
%
\begin{equation}
\label{proof_e_7} \frac{\sqrt{n}}{u_n^{2}k_n} \bigl(H\bigl(p,u_n,\widehat{\beta
}^{fs}\bigr)-H(p,u_n,2) \bigr) = o_p(1).
\end{equation}

We proceed with the rest of the terms in the decomposition of $\widehat
{\mathcal{L}}^n(p,u_n,\beta) - \mathcal{L}(p,u_n,\beta)$ and
$\widehat
{\mathcal{L}}^n(p,v_n,\beta) - \mathcal{L}(p,v_n\beta)$.
We start with the term $R_1^n(u_n)$. It relies on the bound in (\ref
{prelim_a-1}), which in turn depends on the analysis of the term $A_3$
in Section~5.2.3 of \cite{T13}. When $L$ is a Brownian motion, the
bounds for this term get slightly changed. In particular, the bound in
equation (41) of that paper becomes now $K\Delta_n^{1-\iota}$ for
$q>r\vee1$ (this follows by using integration by parts and the
Burkholder--Davis--Gundy inequality) and arbitrarily small $\iota>0$.
Using this, it is easy to show that when $L$ is a Brownian motion, the
bound in (\ref{prelim_a-1}) holds with $\alpha_n$ replaced by $\beta
_n$, where
%
\[
\label{proof_e_8} \beta_n = \frac{\Delta_n^{({3}/{2})(1+(p-1/2)\wedge0-\iota
)}}{\sqrt
{k_n}}\vee\Delta_n^{{1}/{(r\vee1)}-\iota}
\vee\Delta _n^{{p}/{\beta'}\wedge1-{p}/{2}-\iota}\vee\Delta_n^{
{(p+1)}/{(r\vee
1+1)}-\iota}.
\]
Now the bound for $R_1^n(u_n)$ becomes
%
\begin{equation}
\label{proof_e_9} \mathbb{E}\biggl\llvert \frac{R_1^n(u_n)}{nu_n^{2}}\biggr\rrvert \leq K
\biggl( \frac
{\beta_n\vee k_n^{-{1}/{p}+\iota}}{u_n^2} \biggr).
\end{equation}
Further, using the same steps as in the proofs of Lemmas \ref
{lema:prelim-a}--\ref{lema:prelim-c}, as well as
%
\[
\label{proof_e_10} \sup_{u,x\in\mathbb{R}_+} \bigl(|u|^p\bigl|f'_{i,u}(x)\bigr|+|u|^{2p}\bigl|f^{\prime\prime}_{i,u}(x)\bigr|
\bigr)<\infty,
\]
we get
%
\begin{eqnarray}
\label{proof_e_11} &&\mathbb{E}\biggl\llvert \frac{R_2^n(u_n)}{nu_n^{2}}\biggr\rrvert \leq K
\beta _nu_n^{-2},\qquad \mathbb{E}\biggl\llvert
\frac{R_4^n(u_n)}{nu_n^{2}}\biggr\rrvert \leq K(k_n\Delta_n)^{1-\iota},
\\
\label{proof_e_12}
&&\mathbb{E}\biggl\llvert \frac{R_3^n(u_n)}{nu_n^{2}}\biggr\rrvert \leq
Ku_n^{-2-2p} \bigl( (k_n\Delta_n)^{1-\iota}
\vee k_n^{-1/2}(k_n\Delta_n)^{
{1}/{r}\wedge
{(2-p)}/{2}-\iota}
\bigr),
\\
\label{proof_e_13}
&&\frac{\mathbb{E}\llvert \widehat{Z}_1^n(u_n)-\overline
{Z}_1^n(u_n)\rrvert }{nu_n^{2}}\leq K \biggl(
\frac{(\beta_n\vee k_n^{-{1}/{p}+\iota
})}{u_n^{2}} \vee\sqrt{
\Delta_n}(k_n\Delta_n)^{1/2-\iota}
\biggr),
\\
\label{proof_e_14} %
&&\frac{\mathbb{E}\llvert \widehat{Z}_2^n(u_n)-\overline{Z}_2^{(a,n)}(u_n)
- \overline{Z}_2^{(b,n)}(u_n)\rrvert }{nu_n^{2}}
\nonumber
\\[-8pt]
\\[-8pt]
\nonumber
&&\qquad \leq K \biggl(\frac{k_n^{-{1}/{p}+\iota
}}{u_n^{2+2p}} \vee k_n^{-3/2+\iota}\vee
k_n^{-1/2}(k_n\Delta _n)^{
{1}/{r}\wedge{(2-p)}/{2}-\iota}
\biggr). %
\end{eqnarray}
Combining the bounds in (\ref{proof_e_9})--(\ref{proof_e_14}), together
with (\ref{proof_e_3})--(\ref{proof_e_5}), the result in (\ref
{proof_e_6}) and (\ref{proof_e_7}), we establish Lemma~\ref
{lema:prelim-e}. We further note that when $X$ is a L\'{e}vy process,
$R_3^n(u)$ and $R_4^n(u)$ are identically zero.
\end{pf}

We proceed with the proof of Theorem~\ref{thm:cont}. Using Taylor's
expansion and the result in (\ref{prelim_e-1}), $\widehat{Z}_n$,
defined in the statement of Lemma~\ref{lema:prelim-e}, is
asymptotically equivalent to
%
\begin{eqnarray*}
&&\frac{1}{C_{p,2}u_n^2} \bigl(-\log\bigl(\widehat{\mathcal {L}}^n
\bigl(p,u_n,\widehat {\beta}^{fs}\bigr)'\bigr)
- C_{p,2}u_n^2 \bigr) \\
&&\qquad{}- \frac{1}{C_{p,2}v_n^2}
\bigl(-\log\bigl(\widehat{\mathcal{L}}^n\bigl(p,v_n,
\widehat{\beta}^{fs}\bigr)'\bigr) - C_{p,2}v_n^2
\bigr).
\end{eqnarray*}
Using again Taylor's series expansion, the result in (\ref{prelim_e-1})
and that $u_n^{-2}\sqrt{\Delta_n}\rightarrow0$, we have that the above
is asymptotically equivalent to
%
\begin{eqnarray*}
&& \bigl(\log \bigl(-\log\bigl(\widehat{\mathcal{L}}^n
\bigl(p,u_n,\widehat {\beta }^{fs}\bigr)'
\bigr) \bigr)-\log\bigl(C_{p,2}u_n^2\bigr)
\bigr)
\nonumber
\\[-8pt]
\\[-8pt]
\nonumber
&&\qquad - \bigl(\log \bigl(-\log\bigl(\widehat {\mathcal{L}}^n
\bigl(p,v_n,\widehat{\beta}^{fs}\bigr)'\bigr)
\bigr)-\log \bigl(C_{p,2}v_n^2\bigr) \bigr).
\end{eqnarray*}
From here result (\ref{cont_3}) in Theorem~\ref{thm:cont}, both in the
general and L\'{e}vy case, follows from Lemma~\ref{lema:prelim-e}.
\section*{Acknowledgments}
I would
like to thank the Editor, the Associate Editor and two anonymous
referees for many useful suggestions and comments. I would also like to
thank Denis Belomestny, Jose Manuel Corcuera, Valentine Genon-Catalot,
Jean Jacod, Cecilia Mancini, Philip Protter, Markus Reiss, Peter
Spreij, Mathias Vetter and seminar participants at the workshop on
Statistical Inference for L\'{e}vy processes at the Lorentz Center,
University of Leiden and the workshop on Statistics of High-Frequency
Data at Humboldt University.

%
%





\printaddresses
\end{document}